\documentclass[11pt,a4paper,oneside]{article}

\usepackage{hyperref} 
\usepackage{bookmark}

\usepackage{amssymb,amsmath,amsthm}
\usepackage{amsmath}
\usepackage[english]{babel}
\usepackage{graphicx}

\usepackage{wrapfig}
\usepackage{tikz}
	\usetikzlibrary{matrix,arrows}
\usepackage{subfig}
\usepackage{relsize}
\usepackage{xfrac}
\usetikzlibrary{calc}
\usepackage{bbm}
\usepackage{stmaryrd} 

\usepackage{epstopdf}
\usepackage{color}
\usepackage{enumitem}
\usepackage{makeidx}
\makeindex

\theoremstyle{plain}
\newtheorem{theorem}{Theorem}[section]
\newtheorem{lemma}[theorem]{Lemma}

\newtheorem{corollary}[theorem]{Corollary}

\theoremstyle{definition}
\newtheorem{definition}[theorem]{Definition}

\newtheorem{remark}[theorem]{Remark}
\newtheorem{example}[theorem]{Example}
\newtheorem{examples}[theorem]{Examples}
\newtheorem{observation}[theorem]{}

\newcommand{\R}{\mathbb{R}}
\newcommand{\Q}{\mathbb{Q}}
\newcommand{\Z}{\mathbb{Z}}

\newcommand{\N}{\mathbb{N}}
\newcommand{\1}{\mathbbm{1}}

\newcommand{\cA}{\mathcal{A}}
\newcommand{\cB}{\mathcal{B}}

\newcommand{\cI}{\mathcal{I}}

\newcommand{\cL}{\mathcal{L}}
\newcommand{\cM}{\mathcal{M}}
\newcommand{\cN}{\mathcal{N}}

\newcommand{\cP}{\mathcal{P}}

\newcommand{\cX}{\mathcal{X}}

\newcommand{\fB}{\mathfrak{B}} 	\newcommand{\fb}{\mathfrak{b}}

\newcommand{\fP}{\mathfrak{P}} 	\newcommand{\fp}{\mathfrak{p}}

\newcommand{\fU}{\mathfrak{U}}

\newcommand{\sumn}{\sum_{n\in \N}}
\newcommand{\limn}{\lim_{n\rightarrow \infty}}
\newcommand{\limN}{\lim_{N\rightarrow \infty}}
\newcommand{\infn}{\inf_{n\in \N}}
\newcommand{\infm}{\inf_{m\in \N}}
\newcommand{\infk}{\inf_{k\in \N}}
\newcommand{\infN}{\inf_{N\in \N}}

\newcommand{\supn}{\sup_{n\in \N}}

\newcommand{\supk}{\sup_{k\in \N}}
\newcommand{\supN}{\sup_{N\in \N}}

\renewcommand{\epsilon}{\varepsilon}

\DeclareMathOperator{\diam}{diam}
\DeclareMathOperator{\D}{\, \mathrm{d}\!}

\usepackage{geometry}
\usepackage{numname} 

\begin{document}

\renewcommand{\thefootnote}{\Roman{footnote}}

\title{Integrals for functions with values in a partially ordered vector space}

\author{
\renewcommand{\thefootnote}{\Roman{footnote}}
A.C.M. van Rooij
\footnotemark[1]
\\
\renewcommand{\thefootnote}{\Roman{footnote}}
W.B. van Zuijlen
\footnotemark[2]
}

\footnotetext[1]{
Radboud University Nijmegen, Department of Mathematics, P.O.\ Box
9010, 6500 GL Nijmegen, the Netherlands.
}
\footnotetext[2]{
Leiden University, Mathematical Institute, P.O.\ Box 9512, 2300 RA, Leiden, the
Netherlands.
}

\maketitle 

\renewcommand{\thefootnote}{\arabic{footnote}} 

\begin{abstract}
We consider integration of functions with values in a partially ordered vector space, 
and two notions of extension of the space of integrable functions. 
Applying both extensions to the space of real valued simple functions on a measure space leads to the classical space of integrable functions. 

\bigskip\noindent
{\it Key words and phrases.} Partially ordered vector space, Riesz space, Bochner integral, Pettis integral, integral, vertical extension, lateral extension. 
\end{abstract}


\section{Introduction}
\label{section:introduction}

For functions with values in a Banach space there exist several notions of integration.
The best known are the Bochner and Pettis integrals (see \cite{Bo33} and \cite{Pe38}). 
These have been thoroughly studied, yielding a substantial theory (see Chapter III in the book by E. Hille and R.S. Phillips, \cite{HiPh57}). 

As far as we know, there is no notion of integration for functions with values in a partially ordered vector space; not necessarily a $\sigma$-Dedekind complete Riesz space. 
In this paper we present such a notion. 
The basic idea is the following. (Here, $E$ is a partially ordered vector space in which our integrals take their values.)

In the style of Daniell \cite{Da18} and Bourbaki \cite[Chapter 3,4]{Bourbaki52integration},
we do not start from a measure space but from a set $X$, 
a collection $\Gamma$ of functions $X\rightarrow E$, and a functional $\varphi : \Gamma \rightarrow E$, our ``elementary integral''. 
We describe two procedures for extending $\varphi$ to a larger class of functions $X\rightarrow E$. 
The first (see \S\ref{section:vertical_extension}), the ``vertical extension'', is analogous to the usual 
construction of the Riemann integral, proceeding from the space of simple functions. 
The second (see \S\ref{section:lateral_extension}), the ``lateral extension'', is related to the improper Riemann integral. 

In \S\ref{section:combinations_of_extensions} we investigate what happens if one repeatedly applies those extension procedures, without considering the space $E$ to be $\sigma$-Dedekind complete or even Archimedean. 
However, under some mild conditions on $E$ one can embed $E$ into a $\sigma$-Dedekind complete space. 
In \S \ref{section:embeddings_in_larger_space} we discus the extensions procedures in the larger space. 
\S\ref{section:integration_in_R} and \S\ref{section:extensions_of_simple_functions} treat the situation in which $\Gamma$ consists of the simple $E$-valued functions on a measure space. (In \S\ref{section:integration_in_R} we have $E=\R$.)
In \S\ref{section:comparison_B_and_P} we consider connections of our extensions with the Bochner and the Pettis integrals for the case where $E$ is a Banach lattice. 
In \S \ref{section:extensions_of_bochner} we apply our extensions to the Bochner integral.


\section{Some Notation}
\label{section:some_notation}

$\N$ is $\{1,2,3,\dots\}$. \\
Let $X$ be a set. 
We write $\cP(X)$ for the set of subsets of $X$. 
For a subset $A$ of $X$:
\begin{align*}
\label{eqn:indicator_function}
\1_A(x) = \begin{cases}
1 & \mbox{ if } x\in A, \\
0 & \mbox{ if } x\notin A. 
\end{cases}
\end{align*}
As a shorthand notation we write $\1= \1_X$. \\
Let $E$ be a vector space. 
We  write $x= (x_1,x_2,\dots)$ for functions $x: \N \rightarrow E$ (i.e., elements of $E^\N$) and we define
\begin{align*}
c_{00}[E] 
&= \{ x\in E^\N: \exists N \ \forall n\ge N \left[ x_n =0 \right] \}, & c_{00} &= c_{00}[\R]
\end{align*}
We write $c_0$ for the set of sequences in $\R$ that converge to $0$, 
$c$ for the set of convergent sequences in $\R$, $\ell^\infty(X)$ for the set of bounded functions $X\rightarrow \R$, $\ell^\infty$ for $\ell^\infty(\N)$, and $\ell^1$ for the set of absolutely summable sequences in $\R$. 
We write $e_n$ for the element $\1_{\{n\}}$ of $\R^\N$. \\
For a complete $\sigma$-finite measure space $(X,\cA,\mu)$ we write $\cL^1(\mu)$ for the space of integrable functions, $L^1(\mu)= \cL^1(\mu)/\cN$ where $\cN$ denotes the space of functions that are  zero $\mu$-a.e. 
Moreover we write $L^\infty(\mu)$ for the space of equivalence classes of measurable functions that are almost everywhere bounded. 

For a subset $\Gamma$ of a partially ordered vector space $\Omega$, we write $\Gamma^+=\{f\in \Gamma: f \ge 0\}$. 
If $\Lambda, \Upsilon \subset \Omega$ and $f\le g$ for all $f\in \Lambda$ and $g\in \Upsilon$ we write $\Lambda \le \Upsilon$; if $\Lambda =\{f\}$ we write $f \le \Upsilon$ instead of $\{f\} \le \Upsilon$ etc.  
For a sequence $(h_n)_{n\in\N}$ in a partially ordered vector space we write $h_n \downarrow 0$ if $h_1 \ge h_2 \ge h_3 \ge \cdots$ and $\infn h_n=0$.

\section{The vertical extension}
\label{section:vertical_extension}

\emph{\textbf{Throughout this section, $E$ and $\Omega$ are partially ordered vector spaces, $\Gamma \subset \Omega$ is a linear subspace and $\varphi: \Gamma \rightarrow E$ is order preserving and linear. 
Additional assumptions are given in \ref{observation:extra_assumptions_vertical_chapter}.}}

\begin{definition}
\label{phi_integral}
Define 
\begin{align}
\Gamma_v & = \Big \{ f\in \Omega: 
\sup_{\sigma \in \Gamma: \sigma \le f } \varphi(\sigma) =
\inf_{ \tau \in \Gamma: \tau \ge f} \varphi(\tau) 
\Big \},
\end{align}
and $\varphi_v : \Gamma_v \rightarrow E$ by 
\begin{align}
 \varphi_v (f) = \sup_{\sigma \in \Gamma: \sigma \le f } \varphi(\sigma) \qquad (f\in \Gamma_v).
\end{align}
\end{definition}
Note: If $f\in \Omega$ and there exist subsets $\Lambda, \Upsilon \subset \Gamma$ with $\Lambda \le f \le \Upsilon$ such that $ \sup \varphi(\Lambda) = \inf \varphi(\Upsilon)$, then  $f\in \Gamma_v$ and $ \varphi_v (f) =\inf \varphi( \Upsilon) $. 

\begin{observation}
\label{observation:switching_to_subsets}
The following observations are elementary. 
\begin{enumerate}
[label=(\alph*),align=left,leftmargin=1cm,topsep=3pt,itemsep=0pt,itemindent=-1.2em,parsep=0pt]
\item $\Gamma \subset \Gamma_v$ and $\varphi_v(\tau) = \varphi(\tau)$ for all $\tau \in \Gamma$. 
\item $\Gamma_v$ is a partially ordered vector space and $\varphi_v$ is a linear order preserving map\footnote{This follows from the following fact: Let $A,B\subset E$. If $A$ and $B$ have suprema (infima) in $E$, then so does $A+B$ and
$\sup(A+B)= \sup A + \sup B \quad (\inf(A+B)=\inf A + \inf B )$.}. 
\item $(\Gamma_v)_v = \Gamma_v$ and $(\varphi_v)_v= \varphi_v$. 
\item If $\Pi $ is a subset of $\Gamma$, then $\Pi_v \subset \Gamma_v$. 
\end{enumerate}
\end{observation}

Of more importance to us then $\Gamma_v$ and $\varphi_v$ is the following variation in which we consider only countable subsets of $\Gamma$. 

\begin{definition}
\label{def:phi_V_integral}
Let $\Gamma_V$ be the set consisting of those $f$ for which 
there exist countable sets
$\Lambda, \Upsilon \subset \Gamma$ with $\Lambda \le f \le \Upsilon$  such that 
\begin{align}
\label{eqn:inf_and_sup_on_subsets}
 \sup \varphi(\Lambda) = \inf \varphi(\Upsilon).
\end{align}
From the remark following Definition \ref{phi_integral}
it follows that $\Gamma_V$ is a subset of $\Gamma_v$ and that (for $f$ and $\Lambda$ as above) $\varphi_v(f)$ is equal to $\sup \varphi(\Lambda)$.
We will write $\varphi_V=\varphi_{v}|_{\Gamma_V}$. 
We call $\Gamma_V$ the 
\emph{vertical extension}\footnote{One could also define the vertical extension in case $E$, $\Omega$, $\Gamma\subset \Omega$ are  partially ordered sets (not necessarily vector spaces) and $\varphi: \Gamma \rightarrow E$ is an order preserving map. 
}
 under $\varphi$ of $\Gamma$ and $\varphi_V$ the 
\emph{vertical extension} of $\varphi$. 
\end{definition}

In what follows we will only consider $\varphi_V$ and not $\varphi_v$. However, most of the theory presented can be developed similarly for $\varphi_v$. (For comments see \ref{observation:discussion_vertical_extension}.)

\begin{example}
\label{example:Riemann_int_functions_as_vertical_ext}
$\Gamma_V$ is the set of Riemann integrable functions on $[0,1]$ and $\varphi_V$ is the Riemann integral in case $E=\R$ and $\Gamma$ is the linear span of $\{\1_{I}: I \mbox{ is an interval in } [0,1]\}$ and $\varphi$ is the Riemann integral on $\Gamma$. 
\end{example}

\begin{observation}
\label{observation:about_Vertical_ext}
In analogy with \ref{observation:switching_to_subsets} we have the following. 
\begin{enumerate}
[label=(\alph*),align=left,leftmargin=1cm,topsep=3pt,itemsep=0pt,itemindent=-1.2em,parsep=0pt]
\item $\Gamma \subset \Gamma_V $ and $\varphi_V(\tau) = \varphi(\tau)$ for all $\tau \in \Gamma$. 
\item $\Gamma_V$ is a partially ordered vector space and $\varphi_V$ is a linear order preserving map. 
\item $(\Gamma_V)_V = \Gamma_V$ and $(\varphi_V)_V= \varphi_V$. 
\item If $\Pi \subset \Gamma$, then $\Pi_V \subset \Gamma_V$. 
\end{enumerate}
\end{observation}

\begin{definition}
Let $D$ be a linear subspace of $E$. 
$D$ is called 
\emph{mediated} 
 in $E$ if the following is true:
\begin{align}
\label{eqn:mediated}
\notag 
&\mbox{If $A$ and $B$ are countable subsets of $D$ such that $\inf A-B=0$ in $E$, then} \\
&\mbox{$A$ has an infimum (and consequently $B$ has a supremum and $\inf A = \sup B$)}.
\end{align}
$D$ is mediated in $E$ if and only if the following requirement (equivalent with order completeness in the sense of \cite{Pagter81}, for $D=E$) is satisfied 
\begin{align}
\notag &\mbox{If $A$ and $B$ are countable subsets of $D$ such that $\inf A-B=0$ in $E$, then} \\
& \mbox{there exists an $h\in E$ with $B \le h\le A$.}
\label{eqn:order_complete_pagter}
\end{align}
We say that $E$ is \emph{mediated} if $E$ is mediated in itself. 

Note: if $D$ is mediated in $E$, then so is every linear subspace of $D$.
Every $\sigma$-Dedekind complete $E$ is mediated, but so is $\R^2$, ordered lexicographically.
Also, $c_{00}$ and $c_0$ are mediated in $c$, but $c$ is not mediated. 
\end{definition}

With this the following lemma is a tautology. 

\begin{lemma}
\label{lemma:V_extension_in_terms_of_infimum}
Suppose $\varphi(\Gamma)$ is mediated in $E$. Let $f\in \Omega$. 
Then $f\in \Gamma_V$ if and only if there exist countable sets $\Lambda, \Upsilon\subset \Gamma $ with $\Lambda \le f \le \Upsilon$ such that 
\begin{align}
\label{eqn:infimum_differences_in_V_complete_context}
\inf_{\tau \in\Upsilon, \sigma \in \Lambda} \varphi (\tau - \sigma) =0. 
\end{align}
\end{lemma}

The next example shows that $\Gamma_V$ is not necessarily a Riesz space even if $E$ and $\Gamma$ are. However, see Corollary \ref{cor:mediated_E_implies_Gamma_V_Riesz}.

\begin{example}
\label{example:gamma_V_not_Riesz}
Consider $E=c$, $\Gamma= c\times c$, $\Omega = \ell^\infty \times \ell^\infty$. 
Let $\varphi: \Gamma \rightarrow c$ be given by $\varphi (f,g) = f+g$. 
For all $f\in \ell^\infty$ there are $h_1,h_2,\dots \in c$ with $h_n \downarrow f$. 
It follows that, $\Gamma_V = \{(f,g)\in \ell^\infty \times \ell^\infty : f + g \in c\}$. 
Note that $\Gamma_V$ is not a Riesz space since for every $f\in \ell^\infty$ with $f\ge 0$ and $f\notin c$ we have 
$(f,-f)\in \Gamma_V$ but $(f,-f)^+= (f,0)\notin \Gamma_V$. 
\end{example}

\begin{lemma}
\label{lemma:positive_linear_map_maps_Gamma_V_in_itself}
Suppose $\varphi(\Gamma)$ is mediated in $E$.
Let $\Theta: \Omega \rightarrow \Omega$ be an order preserving map with the properties: 
\begin{itemize}[topsep=3pt,itemsep=0pt,parsep=0pt]
\item  if $\sigma, \tau \in \Gamma$ and $ \sigma \le \tau $, then $0 \le \Theta(\tau) - \Theta(\sigma) \le \tau - \sigma$; 
\item $\Theta (\Gamma) \subset \Gamma_V$.   
\end{itemize}
Then $\Theta(\Gamma_V) \subset \Gamma_V$. 
\end{lemma}
\begin{proof}
Let $f\in \Gamma_V$ and 
let $ \Lambda, \Upsilon\subset \Gamma $ be countable sets with $\Lambda \le f \le \Upsilon$ satisfying \eqref{eqn:infimum_differences_in_V_complete_context}. 
Then $\Theta(\Lambda) \le \Theta(f) \le \Theta(\Upsilon)$ and
\begin{align}
\inf_{\tau \in \Theta(\Upsilon), \sigma \in \Theta(\Lambda)} \varphi ( \tau - \sigma)
=
\inf_{\tau \in\Upsilon, \sigma \in \Lambda} \varphi \big( \Theta(\tau) - \Theta( \sigma) \big) 
\le  
\inf_{\tau \in\Upsilon, \sigma \in \Lambda}  \varphi ( \tau- \sigma)
 =0. 
\end{align}
\
\end{proof}

\begin{corollary}
\label{cor:mediated_E_implies_Gamma_V_Riesz}
Suppose that $\varphi(\Gamma)$ is mediated in $E$. 
Suppose $\Omega$ is a Riesz space and $\Gamma$ is a Riesz subspace of $\Omega$. Then so is $\Gamma_V$. 
\end{corollary}
\begin{proof} 
Apply Theorem \ref{lemma:positive_linear_map_maps_Gamma_V_in_itself} with $\Theta(\omega) = \omega^+$.
\end{proof}

\begin{observation}
\label{observation:Gamma_directed_then_Gamma_V_too}
If $\Gamma$ is a directed set, i.e., $\Gamma= \Gamma^+ - \Gamma^+$, then so is $\Gamma_V$. Indeed, if $f\in \Gamma_V$, then there exist $\sigma , \tau\in \Gamma^+$ such that $f \ge \tau - \sigma $ and thus $ f = ( f + \sigma) - \sigma \in \Gamma_V^+ - \Gamma_V^+$.
\end{observation}

\begin{observation}
In the last part of this section we will consider a situation 
in which $\Omega$ has some extra structure. 
But first we briefly consider the case where $E$ is a Banach lattice with $\sigma$-order continuous norm. As it turns out, such an $E$ is mediated (see Theorem \ref{theorem:banach_lattice_sigma_order_c_norm_mediated_and_splitting}), but is not necessarily $\sigma$-Dedekind complete (consider the Banach lattice $C(X)$ where $X$ is the one-point compactification of an uncountable discrete space). 
For such $E$ we describe $\Gamma_V$ in terms of the norm. 
\end{observation}

\begin{theorem}
\label{theorem:vertical_functions_in_banach_lattice_with_sigma_order_cont_norm_context}
Let $E$ be a Banach lattice with a $\sigma$-order continuous norm. Let $\Omega$ be a Riesz space and $\Gamma$ be a Riesz subspace of $\Omega$. 
For $f\in \Omega$ we have: $f\in \Gamma_V$ if and only if for every $\epsilon>0$ there exist $\sigma, \tau \in \Gamma$ with $\sigma \le f \le \tau$ and $\|\varphi(\tau) - \varphi(\sigma) \|<\epsilon$. 
\end{theorem}
\begin{proof}
First, assume $f\in \Gamma_V$. 
As $\Gamma$ is a Riesz subspace of $\Omega$ there exist sequences $(\sigma_n)_{n\in\N}$ and $(\tau_n)_{n\in\N}$ in $\Gamma$ such that $\sigma_n \uparrow$, $\tau_n \downarrow$, 
\begin{align}
\sigma_n \le f \le \tau_n \qquad (n\in\N), \qquad \supn \varphi(\sigma_n )= \infn \varphi(\tau_n). 
\end{align}
Then $\varphi(\tau_n - \sigma_n) \downarrow 0$ in $E$, so $\|\varphi(\tau_n)-\varphi(\sigma_n)\|\downarrow 0$ and we are done. \\
The converse: For each $n\in\N$, choose $\sigma_n,\tau_n\in \Gamma$ for which 
\begin{align}
\sigma_n \le f \le \tau_n, \qquad \|\varphi(\tau_n) - \varphi(\sigma_n)\| \le n^{-1}. 
\end{align}
Setting $\sigma_n' = \sigma_1 \vee \cdots \vee \sigma_n$ and $\tau_n'= \tau_1 \wedge \cdots \wedge \tau_n$ we have, for each $n\in\N$
\begin{align}
\sigma_n', \tau_n'\in \Gamma, \qquad \sigma_n' \le f \le \tau_n'.
\end{align}
If $n\ge N$, then $0\le \sigma_n' - \sigma_N' \le f - \sigma_N \le \tau_N - \sigma_N$, whence $\|\varphi(\sigma_n') - \varphi(\sigma_N')\| \le \|\varphi(\tau_N) - \varphi(\sigma_N) \| \le N^{-1}$. 
Thus, the sequence $(\varphi(\sigma_n'))_{n\in\N}$ converges in the sense of the norm. 
So does $(\varphi(\tau_n'))_{n\in\N}$. 
Their limits are the same element $a$ of $E$, and, since $\sigma_n' \uparrow, \tau_m' \downarrow$, we see that $a= \supn \varphi(\sigma_n') = \infm \varphi(\tau_m')$.
\end{proof}

\begin{observation}
\label{observation:extra_assumptions_vertical_chapter}
\textbf{\emph{In the rest of this section $\Omega$ is the collection $F^X$  of all maps of a set $X$ into a partially ordered vector space $F$.}} 
\end{observation}

\begin{observation}
A function $g: X \rightarrow \R$ determines a multiplication operator $f \mapsto gf$ in $\Omega$. We investigate the collection of all functions $g$ for which 
\begin{align}
f\in \Gamma_V \quad \Longrightarrow \quad gf \in \Gamma_V,
\end{align}
and, for given $f$, the behaviour of the map $g \mapsto \varphi_V(gf)$. 
\end{observation}

\begin{observation}
\label{observation:order_limits_step_functions}
For an algebra of subsets of $X$, $\cA\subset \cP(X)$ we write
$[ \cA ]$ for the Riesz space of all $\cA$-step functions, i.e., functions of the form $\sum_{i=1}^n \lambda_i \1_{A_i}$ for $n\in\N$, $\lambda_i\in \R$, $A_i \in \cA$ for $i\in\{1,\dots,n\}$. 
Define the collection of functions $[\cA]^o$ by 
\begin{align}
[\cA]^o = & \{ f\in \R^X: \mbox{ there are } (s_n)_{n\in\N} \mbox{ in } [\cA] \mbox{ and } (j_n)_{n\in\N} \mbox{ in } [\cA]^+ \\
\notag &  \mbox{ for which } |f-s_n|\le j_n \mbox{ and } j_n \downarrow 0 \mbox{ pointwise}\}. 
\end{align}
(This $[\cA]^o$ is the vertical extension of $[\cA]$ obtained by, in Definition \ref{def:phi_V_integral}, choosing $E=\R^X,\Omega= \R^X, \Gamma = [\cA], \ \varphi(f)=f \quad (f\in \Gamma)$.)
Note that $[\cA]$ and $[\cA]^o$ are Riesz spaces, and uniform limits of elements of $[\cA]$ are in $[\cA]^o$. 
(Actually, $[\cA]^o$ is uniformly complete.)
Furthermore, $[\cA]^o$ contains every bounded function $f$ with $\{x\in X : f(x) \le s\} \in \cA$ for all $s\in \R$. 
In case $\cA$ is a $\sigma$-algebra, $[\cA]^o$ is precisely the collection of all bounded $\cA$-measurable functions. 
\end{observation}

\begin{lemma}
\label{lemma:seq_of_order_step_f_with_inf_zero}
Let $\cA\subset \cP(X)$ be an algebra of subsets of a set $X$.
Suppose that $(g_n)_{n\in\N}$ is a sequence in $[\cA]^o$ for which $g_n \downarrow 0$ pointwise. Then there exists a sequence $(j_n)_{n\in\N}$ in $[\cA]$ with $j_n \ge g_n$ and $j_n \downarrow 0$ pointwise. 
\end{lemma}
\begin{proof}
For all $n\in\N$ there exists a sequence $(s_{nk})_{k\in\N}$ in $[\cA]$ with $s_{nk} \ge g_n$ for all $k\in\N$ and  $s_{nk} \downarrow_k g_n$ pointwise. 
Since $(g_n)_{n\in\N}$ is a decreasing sequence, we have $s_{mk} \ge g_n$ for all $m\le n$ and all $k\in\N$. Hence $j_n:=\inf_{m,k\le n } s_{mk}$ is an element in $[\cA]$ with $j_n \ge g_n$. Clearly $j_n \downarrow$ and $\infn j_n = \infn \inf_{m,k\le n} s_{mk} = \infn \infk s_{nk} = \infn g_n =0$. 
\end{proof}

The following lemma is a consequence of Theorem \ref{lemma:positive_linear_map_maps_Gamma_V_in_itself}. 

\begin{lemma}
\label{lemma:bounded_measurable_functions}
Define the algebra 
\begin{align}
\label{eqn:cA}
\cA= \{ A \subset X: f\1_A \in \Gamma \mbox{ for } f\in \Gamma\}. 
\end{align}
If $\varphi(\Gamma)$ is mediated in $E$, then
\begin{align}
 f\1_A \in \Gamma_V \qquad (f\in \Gamma_V, A\in \cA). 
\end{align}
\end{lemma}

\begin{definition}
\label{definition:integrally_closed}
$E$ is called \emph{integrally closed} (see Birkhoff \cite{Bi67}) if for all $a,b\in E$ the following holds: if $na \le b$ for all $n\in\N$, then $a\le 0$. 
\end{definition}

\begin{definition}
A sequence $(a_n)_{n\in\N}$ in $E$ is called 
\emph{order convergent}
to an element $a\in E$ if there exists a sequence $(h_n)_{n\in\N}$ in $E^+$ with $h_n \downarrow 0$ and $- h_n \le a-a_n\le h_n$. \\
Notation: $a_n \xrightarrow{o} a$. 
\end{definition}

\begin{theorem}
\label{theorem:bounded_measurable_functions}
Let $\cA$ be as in \eqref{eqn:cA}.
Suppose that $E$ is integrally closed, $\Gamma$ is directed and $\varphi(\Gamma)$ is mediated in $E$. 
Furthermore assume $\varphi$ has the following continuity property. 
\begin{align}
\label{eqn:continuity_property_of_varphi}
&\mbox{If $A_1,A_2,\dots$ in $\cA$ are such that $A_1 \supset A_2 \supset \cdots$ and $\bigcap_{n\in\N} A_n = \emptyset$,} \\
\notag &\mbox{then $\varphi(f\1_{A_n}) \downarrow 0$ for all $f\in \Gamma^+$.
}
\end{align}
\begin{enumerate}
[label=\emph{(\alph*)},align=left,leftmargin=1cm,topsep=3pt,itemsep=0pt,itemindent=-1.2em,parsep=0pt]
\item 
\label{item:gf_in_Gamma_V}
$gf\in \Gamma_V$ for all $g\in [\cA]^o$ and all $f\in \Gamma_V$. 
\item 
\label{item:gf_by_order_convergence}
Let $g\in [\cA]^o$ and let $(g_n)_{n\in\N}$ be a sequence in $[\cA]^o$ for which  
there is a sequence $(j_n)_{n\in\N}$ in $[\cA]^{o+}$ with $-j_n \le g_n -g \le j_n$ and $j_n \downarrow 0$ pointwise. 
Then 
\end{enumerate}
\begin{align}
\label{eqn:order_convergence_bounded_measurable_times_f}
 \varphi_V (g_n f) \xrightarrow{o} \varphi_V(gf) \qquad (f\in \Gamma_V).
\end{align}
(Order convergence in the sense of $E$.)
\end{theorem}
\begin{proof}
We first prove the following: \\
($\star$)
Let $f\in \Gamma_V^+$. 
Let $(g_n)_{n\in\N}$ be a sequence in $[\mathcal{A}]^o$ for which $g_n f \in \Gamma_V$ for all $n\in\N$ and $g_n \downarrow 0$ pointwise. Then 
\begin{align}
 \varphi_V (g_n f) \downarrow 0.
\end{align} 
Let $\sigma \in \Gamma^+$, $\sigma \ge f$. 
It follows from Lemma \ref{lemma:seq_of_order_step_f_with_inf_zero} that we may assume $g_n \in [\cA]$ for all $n\in\N$. 
For all $n\in\N$ we have $0 \le \varphi_V(g_n f) \le \varphi_V( g_n \sigma )$, so we are done if $\varphi_V(g_n \sigma) \downarrow 0$. \\
Let $h\in E$, $h\le \varphi_V(g_n \sigma)$ for all $n\in\N$; we prove $h\le 0$. \\
Take $\epsilon>0$. For each $n\in\N$, set $A_n = \{x\in X: g_n(x) \ge \epsilon\}$. Then $A_n\in \cA$ for $n\in\N$ and $A_1 \supset A_2 \supset \cdots$ and $\bigcap_{n\in\N} A_n = \emptyset$. Putting $M= \|g_1\|_\infty$ we see that 
\begin{align}
 g_n \le \epsilon \1_X + M \1_{A_n} \qquad (n\in\N),
\end{align}
whence 
\begin{align}
 h\le \varphi_V(g_n \sigma) \le \epsilon \varphi( \sigma ) + M \varphi(\1_{A_n} \sigma) \qquad (n\in\N). 
\end{align}
By the continuity property of $\varphi$, $h\le \epsilon \varphi(\sigma)$. As this is true for each $\epsilon>0$ and $E$ is integrally closed, we obtain $h\le 0$. 

\emph{\ref{item:gf_in_Gamma_V}} Since $\Gamma_V$ is directed (see \ref{observation:Gamma_directed_then_Gamma_V_too}) it is sufficient to consider $f\in \Gamma_V^+$.  
Let  $g\in [\cA]^o$. There are sequences of step functions $(h_n)_{n\in\N}$ and $(j_n)_{n\in\N}$ for which $h_n \uparrow g$, $j_n \downarrow g$ and thus $j_n - h_n \downarrow 0$.
By Lemma \ref{lemma:bounded_measurable_functions} $h_n f, j_n f \in \Gamma_V$ for all $n\in\N$. 
Then $h_n f \le g f \le j_n f$ for $n\in\N$ and $\inf_{n\in\N}  \varphi_V((j_n - h_n)f)=0$ by ($\star$). 
By Lemma \ref{lemma:V_extension_in_terms_of_infimum} and \ref{observation:about_Vertical_ext}(c) we obtain that $gf\in \Gamma_V$. 

\emph{\ref{item:gf_by_order_convergence}} It is sufficient to consider $f\in \Gamma_V^+$. By \emph{\ref{item:gf_in_Gamma_V}} we may also assume $g=0$. But then \emph{\ref{item:gf_by_order_convergence}} follows from ($\star$). 
\end{proof}

\begin{remark}
\label{remark:bounded_measurable_are_subset_of_A_V}
Consider the situation in Theorem \ref{theorem:bounded_measurable_functions}. 
Suppose $\cB\subset \cA$ is a $\sigma$-algebra.
Then all bounded $\cB$-measurable functions lie in $[\cA]^o$. 
If  $(g_n)_{n\in\N}$ is a bounded sequence of bounded $\cB$-measurable functions 
that converges pointwise to a function $g$, then the condition of Theorem \ref{theorem:bounded_measurable_functions}(b) is satisfied. 
\end{remark}

\begin{remark}
\label{remark:situations_compared_with_next_chapter}
In the next section we will consider a situation similar to the one of Theorem \ref{theorem:bounded_measurable_functions}, in which $\cA$ is replaced by a subset $\cI$ that is closed under taking finite intersections. We will also adapt the continuity property on $\varphi$ (see \ref{observation:comparison_situation_theorem_and_laterally_extendable}). 
\end{remark}

\section{The lateral extension}
\label{section:lateral_extension}

The construction described in Definition \ref{def:phi_V_integral} is reminiscent of the Riemann integral and, indeed, the Riemann integral is a special case (see Example \ref{example:Riemann_int_functions_as_vertical_ext}). 

In the present section we consider a type of extension, analogous to the \emph{improper} Riemann integral. One usually defines the improper integral of a function $f$ on $[0,\infty)$ to be 
\begin{align}
\lim_{s\rightarrow \infty} \int_0^s f(x) \D x,
\end{align}
approximating the domain, not the values of $f$. 

For our purposes a more convenient description of the same integral would be
\begin{align}
\sum_{n=1}^\infty \int_{a_n}^{a_{n+1}} f(x) \D x,
\end{align}
where $0=a_1<a_2<\cdots$ and $a_n \rightarrow \infty$. Here the domain is split up into manageable pieces. 

Splitting up the domain is the basic idea we develop in this section. (This may explain our use of the terms ``vertical'' and ``lateral''.)

\textbf{\emph{Throughout this section, $E$ and $F$ are partially ordered vector spaces, $\Gamma$ is a directed\footnote{For the construction of the lateral extension, one does not need to assume that $\Gamma$ is directed. However, as one can see later on in  the construction, the only part of $\Gamma$ that matters for the extension is $\Gamma^+ - \Gamma^+$.}
 linear subspace of $F^X$, and $\varphi$ is a linear order preserving map $\Gamma \rightarrow E$. (With $\Omega = F^X$,  all considerations of \S\ref{section:vertical_extension} are applicable.)}}

\begin{center}
\begin{tikzpicture}
\matrix(a)[matrix of math nodes, row sep=1.4em, column sep=-0.5em,
text height=1.5ex, text depth=0.25ex]
{ \Gamma 	&  \subset F^X   \\
E	& \\};
\path[->,font=\scriptsize] (a-1-1) edge node[left]{$\varphi$} (a-2-1);
\end{tikzpicture} 
\end{center}

\textbf{\emph{
Furthermore, $\cI$ is a collection of subsets of $X$, closed under taking finite intersections. See Definition \ref{def:countable_partition} and Definition \ref{def:laterally_extendable} for two more assumptions. 
}}

As a shorthand notation, if $(a_n)_{n\in\N}$ is a sequence in $E^+$ and $\{ \sum_{n=1}^N a_n : N\in\N\}$ has a supremum, we denote this supremum by 
\begin{align}
\label{eqn:sum_n_a_n}
\sum_n a_n.
\end{align}

\begin{definition}
\label{def:countable_partition}
A disjoint sequence $(A_n)_{n\in\N}$ 
of elements in $\cI$ whose union is $X$ is called a 
\emph{partition}\index{partition}. 
If $(A_n)_{n\in\N}$ and $(B_n)_{n\in\N}$ are partitions and for all $n\in\N$ there exists an $m\in\N$ for which $B_n \subset A_m$, then $(B_n)_{n\in\N}$ is called a 
\emph{refinement}\index{refinement}
of $(A_n)_{n\in\N}$. 
Note that if $(A_n)_{n\in\N}$ and $(B_n)_{n\in\N}$ are partitions then there exists a refinement of both $(A_n)_{n\in\N}$ and $(B_n)_{n\in\N}$ (e.g., a partition that consists of all sets of the form $A_n \cap B_m$ with $n,m\in\N$).\\
\textbf{\emph{We assume that there exists at least one partition.}}
\end{definition}

\begin{definition}
\label{def:laterally_extendable}
We call a linear subspace $\Delta$ of $F^X$  
\emph{stable}\index{stable} (under $\cI$)
if 
\begin{align}
\label{eqn:lateral_assumption_1}
 f\1_A\in \Delta \qquad (f\in \Delta, A\in \cI).
\end{align} 
If $\Delta$ is a stable space, then a linear and order preserving map $\omega : \Delta \rightarrow E$ is said to be 
\emph{laterally extendable}\index{laterally extendable}
if for all partitions $(A_n)_{n\in\N}$ 
\begin{align}
\label{eqn:lateral_assumption_2}
 \omega(f) = \sum_n \omega(f\1_{A_n}) 
\quad  \left( \mbox{see} \eqref{eqn:sum_n_a_n} \right) 
 \qquad (f\in \Delta^+).
\end{align} 
\textbf{\emph{We assume $\Gamma$ is stable and $\varphi$ is laterally extendable.}}
\end{definition}

\begin{observation}
\label{observation:comparison_situation_theorem_and_laterally_extendable}
In the situation of Theorem \ref{theorem:bounded_measurable_functions} we can choose $\cI=\cA$; then \eqref{eqn:continuity_property_of_varphi} is precisely the lateral extendability of $\varphi$.
\end{observation}

\begin{example}
\label{example:example_latterally_extendable_sequences}
For any partially ordered vector space $F$ and a linear subspace $E\subset F$, 
the following 
choices lead to a system fulfilling all of our assumptions: 
$X= \N$, $\cI = \cP(\N)$, $\Gamma = c_{00}[E]$ (see \S\ref{section:some_notation}), $\varphi(f) = \sum_{n\in\N} f(n)$  for $f\in \Gamma$. 
\end{example}

\begin{definition}
\label{def:sigma_extension_function}
Let $\Delta$ be a stable subspace of $F^X$ and let $\omega: \Delta \rightarrow E$ be a laterally extendable linear order preserving map.
Let $(A_n)_{n\in\N}$ be a partition, and $f: X \rightarrow F$. 
We call $(A_n)_{n\in\N}$ a 
\emph{partition for} $f$ (occasionally \emph{$\Delta$-partition for} $f$) if 
\begin{align}
\label{eqn:f_on_partition_in_Gamma}
 f \1_{A_n} \in \Delta \qquad (n\in\N).
\end{align}
A function $f:X \rightarrow F$ is said to be a
\emph{partially in}\index{partially in}
$\Delta$ if there exists a partition for $f$. 
For $f: X \rightarrow F^+$, 
$(A_n)_{n\in\N}$ is called a 
\emph{$\omega$-partition for} $f$ 
if it is a partition for $f$ and if 
\begin{align}
\label{eqn:partial_sums_over_partition_has_sup}
 \sum_{n} \omega(f \1_{A_n}) \mbox{ exists}.
\end{align}
A function $f: X \rightarrow F^+$ that is partially in $\Delta$ is called 
 \emph{laterally $\omega$-integrable}
if there exists a $\omega$-partition for $f$.
\end{definition}

\begin{example}
\label{example:example_latterally_extendable_sequences_explained}
Consider the situation of Example \ref{example:example_latterally_extendable_sequences}.
A function $x: \N \rightarrow F$ is partially in $\Gamma$ if and only if $x_n \in E$ for every $n\in\N$. 
If $x \ge 0$, then $x$ is laterally integrable if $x_n \in E$ for every $n\in\N$ and $\sum_n x_n$ exists in $E$.
\end{example}

\begin{observation}
Naturally, we wish to use \eqref{eqn:partial_sums_over_partition_has_sup} to define an integral for $f$. For that we have to show the supremum to be independent of the choice of the partition $(A_n)_{n\in\N}$. 
\end{observation}

\begin{lemma}
\label{lemma:relations_for_partitions}
\
\begin{enumerate}[label=\emph{(\alph*)},align=left,leftmargin=1cm,topsep=3pt,itemsep=0pt,itemindent=-1.2em,parsep=0pt]
\item Let $f: X \rightarrow F$ and let 
$(A_n)_{n\in\N} $ be a partition for $f$. 
If $(B_n)_{n\in\N}$ is a partition that is a refinement of $(A_n)_{n\in\N}$, then $ (B_n)_{n\in\N} $ is a partition for $f$.
\item Let $f: X \rightarrow F^+$ and let $(A_n)_{n\in\N} $ and $ (B_m)_{m\in\N} $ be partitions for $f$. Then the sets
\end{enumerate}
\begin{align}
\left\{ \sum_{n=1}^N \varphi(f\1_{A_n}) : N\in\N \right\} 
\mbox{ and } 
\left\{ \sum_{m=1}^M \varphi(f\1_{B_m}) : M\in\N \right\}
\end{align}
have the same upper bounds in $E$. 
\end{lemma}
\begin{proof}
We leave the proof of (a) to the reader. 
Let $u$ be an upper bound for the set $\{ \sum_{n=1}^N \varphi(f\1_{A_n}) : N\in\N \}$; it suffices to prove that $u$ is an upper bound for $\{\sum_{m=1}^M \varphi(f\1_{B_m}) : M\in\N \}$. 
Take $M\in\N$; we are done if $u \ge \sum_{m=1}^M \varphi(f\1_{B_m})$, i.e., if $u \ge \varphi(f\1_B)$ where $B= B_1 \cup \cdots \cup B_M$. 
But $f\1_B\in \Gamma$ so $\varphi(f\1_B) = \sum_n \varphi(f\1_B \1_{A_n}) = \supN \sum_{n=1}^N \varphi(f\1_B \1_{A_n})$, whereas, for each $N\in\N$ 
\begin{align}
\sum_{n=1}^N \varphi(f\1_B \1_{A_n}) \le \sum_{n=1}^N \varphi(f\1_{A_n}) \le u. 
\end{align}
\
\end{proof}

\begin{theorem}
\label{theorem:unique_integral_elements_of_gamma_sigma}
Let $f:X \rightarrow F^+$ be laterally $\varphi$-integrable. 
Then every partition for $f$ is a $\varphi$-partition for $f$. 
There exists an $a\in E^+$ such that for every partition $(A_n)_{n\in\N}$ for $f$,
 \begin{align}
 \label{eqn:welldefinedness_integral_extention}
   a=  \sum_n \varphi(f \1_{A_n}).  
 \end{align}
 If $f\in \Gamma^+$, then $a= \varphi(f)$. 
\end{theorem}
\begin{proof}
This is a consequence of Lemma \ref{lemma:relations_for_partitions}(b).
\end{proof}

\begin{definition}
\label{def:definition_of_varphi_L_integral_for_positive_function}
For a laterally $\varphi$-integrable $f: X \rightarrow F^+$ we call the element $a\in E^+$ for which \eqref{eqn:welldefinedness_integral_extention} holds its \emph{$\varphi_L$-integral} and denote it by $ \varphi_L(f)$.
For the moment, denote by $(\Gamma^+)_L$ the set of all laterally $\varphi$-integrable functions $f: X \rightarrow F^+$. We proceed to extend $\varphi_L$ to a linear  function defined on the linear hull of $(\Gamma^+)_L$, see Definition \ref{def:lat_varphi-int}. \end{definition}

\begin{observation}
\label{observation:lateral_integral_well_defined}
The assumptions that $\Gamma$ is stable and $\varphi$ is laterally extendable are crucial for the fact that the $\varphi_L$-integral of a laterally $\varphi$-integrable function is independent of the choice of a $\varphi$-partition (see Lemma \ref{lemma:relations_for_partitions}(b)). 
\end{observation}

\begin{observation}
\label{observation:two_supremum_properties}
We will use the following rules for a partially ordered vector space $E$: 
\begin{align}
\label{eqn:sup_of_sum_is_sum_of_sup}
 a_n \uparrow a, b_n \uparrow b & \Longrightarrow a_n + b_n \uparrow a + b & (a_n, b_n,a,b \in E),  \\ 
\label{eqn:subtraction_Vule_with_sums_of_suprema}
 a_n \uparrow, \ b_n \uparrow b, \ a_n + b_n \uparrow a + b &\Longrightarrow a_n \uparrow a & (a_n, b_n,a,b \in E). 
\end{align}
\end{observation}

\begin{observation}\textbf{(Extending $\varphi_L$)} 
\label{observation:extending_varphi_L}
Define $\Gamma_L= \{ f_1-f_2 : f_1,f_2\in (\Gamma^+)_L\}$. \\
\emph{Step 1.} 
Let $f,g\in (\Gamma^+)_L$. There exists an $(A_n)_{n\in\N}$ that is a $\varphi$-partition for $f$ and for $g$. By defining $a_N= \sum_{n=1}^N \varphi(f\1_{A_n})$ and $b_N=\sum_{n=1}^N \varphi(g\1_{A_n})$ for $N\in\N$, by \eqref{eqn:sup_of_sum_is_sum_of_sup} we obtain $f+g\in (\Gamma^+)_L$ with $\varphi_L(f+g) = \varphi_L(f) + \varphi_L(g)$.

 Consequently, $\Gamma_L$ is a vector space, containing $(\Gamma^+)_L$. \\
\emph{Step 2.}
If $g_1,g_2,h_1,h_2\in (\Gamma^+)_L$ and $g_1 -g_2 = h_1 - h_2$, then $g_1 + h_2 = g_2 + h_1$ so that, by the above, $\varphi_L(g_1) - \varphi_L(h_1) = \varphi_L(g_2) - \varphi_L(h_2)$.

Hence, $\varphi_L$ extends to a linear function $\Gamma_L \rightarrow E$ (also denoted by $\varphi_L$). \\
\emph{Step 3.}
Let $f,g \in (\Gamma^+)_L$ and $f\le g$. By defining $a_N$ and $b_N$ as in step 1 and $c_N=b_N - a_N$, by \eqref{eqn:subtraction_Vule_with_sums_of_suprema} we infer that $g-f\in (\Gamma^+)_L$. 

Thus, if $f\in \Gamma_L$ and $f\ge 0$, then $f\in (\Gamma^+)_L$. Briefly: $(\Gamma^+)_L$ is $ \Gamma_L^+$, the positive part of $\Gamma_L$. 
\end{observation}

\begin{definition}
\label{def:lat_varphi-int}
A function $f: X \rightarrow F$ is called 
\emph{laterally $\varphi$-integrable}
 if $f\in \Gamma_L$
  (see \ref{observation:extending_varphi_L}), i.e., if 
 there exist $f_1,f_2\in (\Gamma^+)_L$ for which $f= f_1 - f_2$. 
The $\varphi_L$-integral of such a function is defined by 
$\varphi_L(f) = \varphi_L(f_1) - \varphi_L(f_2)$. \\
$\varphi_L$ is a function $ \Gamma_L \rightarrow E$ and is called the 
\emph{lateral extension} of $\varphi$. 
The set of laterally $\varphi$-integrable functions, $\Gamma_L$, is called the \emph{lateral extension} of $\Gamma$ under $\varphi$. \\
Note that, thanks to Step 3 of \ref{observation:extending_varphi_L}, this definition of ``laterally $\varphi$-integrable'' does not conflict with the one given in Definition \ref{def:definition_of_varphi_L_integral_for_positive_function}. 
\end{definition}

\begin{observation}
\label{observation:some_observations_regarding_the_lateral_extension}
Like for the vertical extension, we have the following elementary observations:
\begin{enumerate}
[label=(\alph*),align=left,leftmargin=1cm,topsep=3pt,itemsep=0pt,itemindent=-1.2em,parsep=0pt]
\item $\Gamma \subset \Gamma_L$\footnote{Note that for this inclusion it is necessary that $\Gamma$ be directed.} and $\varphi_L(\tau) = \varphi(\tau)$ for all $\tau \in \Gamma$. 
\item $\Gamma_L$ is a directed partially ordered vector space and $ \varphi_L$ is a linear order preserving function on $\Gamma_L$. 
\item If $\Pi$ is a directed linear subspace of $F^X$ and $\Pi \subset \Gamma$, then $\Pi_L \subset \Gamma_L$. 
\end{enumerate}
($(\Gamma_L)_L$ is not so easy. See Theorem \ref{theorem:stability_implies_laterally_extendability_for_V_and_L} and Example \ref{example:gamma_L_not_stable}.)
\end{observation}

In case $E$ is a Banach lattice with $\sigma$-order continuous norm, for $\Gamma_L^+$ we have an analogue of Theorem \ref{theorem:vertical_functions_in_banach_lattice_with_sigma_order_cont_norm_context}. 

\begin{lemma}
\label{lemma:lateral_extended_functions_in_banach_lattice_with_sigma_order_cont_norm_context}
Suppose $E$ is a Banach lattice with $\sigma$-order continuous norm. 
Let $f: X \rightarrow F^+$. 
Then $f$ lies in $\Gamma_L^+$ if and only if there exists a $\Gamma$-partition $(A_n)_{n\in\N}$ for $f$ such that the sequence $(\varphi (f\1_{A_n}))_{n\in\N}$ has a sum in the sense of the norm, in which case $\varphi_L(f)$ is this sum. 
\end{lemma}
\begin{proof}
The ``only if'' part follows by definition of $\Gamma_L$ and the $\sigma$-order continuity of the norm. 
For the ``if'' part; this follows from the fact that if $a_n \uparrow$ and $\|a_n - a\|\rightarrow 0$ for $a,a_1,a_2,\dots \in E$, then $a_n \uparrow a$. 
\end{proof}

We will now investigate conditions under which $\varphi_L$ and $\varphi_V$ themselves are laterally extendable.
(For that, their domains have to be able to play the role of $\Gamma$, so they have to be stable.)
First a useful lemma: 

\begin{lemma}
\label{lemma:a_tool_for_refinements}
Let $f\in \Gamma_L$. 
Then there exists a partition $(A_n)_{n\in\N}$ for $f$ such that 
every refinement $(B_m)_{m\in\N}$ of it (is a partition for $f$ and) has this property:
\begin{align}
\label{eqn:tool_refinements_property}
h\in E, \  h \ge \sum_{m=1}^M \varphi(f\1_{B_m}) \mbox{ for all } M\in \N 
\quad \Longrightarrow \quad 
h\ge \varphi_L(f). 
\end{align}
\end{lemma}
\begin{proof}
Write $f= f_1 - f_2$ with $f_1, f_2 \in \Gamma_L^+$. Let $(A_n)_{n\in\N}$ be a partition for $f_1$ and $f_2$, and let $(B_m)_{m\in\N}$ be a refinement of $(A_n)_{n\in\N}$. 
Note that $(B_m)_{m\in\N}$ is a partition for $f_1$ and $f_2$. 
Let $h$ be an upper bound for $\{ \sum_{m=1}^M \varphi(f\1_{B_m}): M\in\N\}$ in $E$. 
For all $M\in\N$, 
\begin{align}
h+ \sum_{m=1}^M \varphi (f_2 \1_{B_m}) \ge \sum_{m=1}^M \varphi(f\1_{B_m}) + \sum_{m=1}^M \varphi(f_2 \1_{B_m}) = \sum_{m=1}^M \varphi(f_1 \1_{B_m}). 
\end{align}
Taking the supremum over $M$ yields $h + \varphi_L(f_2) \ge \varphi_L(f_1)$, i.e., $h\ge \varphi_L(f)$. 
\end{proof}

\begin{theorem} \mbox{} 
\label{theorem:stability_implies_laterally_extendability_for_V_and_L}
\begin{enumerate}[label=
\emph{(\alph*)},align=left,leftmargin=1cm,topsep=3pt,itemsep=0pt,itemindent=-1.2em,parsep=0pt]
\item 
Suppose $\Gamma_L$ is stable. 
Then $\varphi_L$ is laterally extendable,  i.e., 
\begin{align}
\varphi_L(f) = \sum_n \varphi_L (f \1_{A_n})
\end{align}
for every $f\in \Gamma_L^+$ and every $\varphi_L$-partition $(A_n)_{n\in\N}$ for $f$. 
Therefore $\left(\Gamma_L\right)_L= \Gamma_L$ and $(\varphi_L)_L = \varphi_L$. 
\item 
Suppose $\Gamma_V$ is stable.
Then $\varphi_V$ is laterally extendable. 
(For $(\Gamma_V)_L$ see \S\ref{section:combinations_of_extensions}.)
\end{enumerate}
\end{theorem}
\begin{proof}
(a)
Let $f\in \Gamma_L^+$ and let $(B_n)_{n\in\N}$ be a $\varphi_L$-partition for $f$. 
Let $(A_n)_{n\in\N}$ be the partition for $f$ as in Lemma \ref{lemma:a_tool_for_refinements}. 
Then form a common refinement of $(B_n)_{n\in\N}$ and $(A_n)_{n\in\N}$ and apply Lemma \ref{lemma:a_tool_for_refinements}. 

(b) Let $f\in \Gamma_V^+$ and let $(A_n)_{n\in\N}$ be a partition. 
Let $h\in E, h \ge \sum_{n=1}^N \varphi_V(f \1_{A_n})$ for every $N\in\N$. 
We wish to prove $h \ge \varphi_V(f)$, which will be the case if $h \ge \varphi(\sigma)$ for every $\sigma \in \Gamma$ with $\sigma \le f$. For that apply Lemma \ref{lemma:a_tool_for_refinements}
 to $\sigma$. 
\end{proof}

The following shows that $\Gamma_L$ may not be stable, in which case there is no $(\Gamma_L)_L$. 
(However, see Theorem \ref{theorem:varphi_of_Gamma_R-closed_then_L_and_V_stable}(a).) 

\begin{example}
\label{example:gamma_L_not_stable}
Consider the situation in Example \ref{example:example_latterally_extendable_sequences} and assume 
there is an $a: \N \rightarrow E^+$ such that $\sum_n a_n$ exists in $F$ and $\sum_n a_{2n}$ does not (e.g. $E=F=c$ and $a_n=e_n = \1_{\{n\}}$). 
By Example \ref{example:example_latterally_extendable_sequences_explained} $a$ lies in $\Gamma_L$ but $b= (0,a_2,0,a_4,\dots)$ does not; but $b = a \1_{\{2,4,6,\dots\}}$ and $\{2,4,6,\dots\}\in \cI$. 
(Actually, the existence of such an $a: \N \rightarrow E^+$ is equivalent to $E$ not being ``splitting'' in $F$; see Definition \ref{def:splitting} and \eqref{eqn:splitting2}.)
\end{example}

\begin{remark}
$\Gamma_V$ may not be stable either. 
With $E=c$, $F= \ell^\infty$, $X= \{1,2\}$, $\Gamma= c\times c$ and $\varphi(f,g)=f +g$ 
(as in Example \ref{example:gamma_V_not_Riesz}), the space $\Gamma_V$ is not stable for $\cI = \cP(X)$.  
\end{remark}

\begin{definition}
\label{def:splitting}
Let $D$ be a linear subspace of $E$. $D$ is called \emph{splitting} in $E$ if the following is true:
\begin{align}
\label{eqn:splitting}
\notag 
&\mbox{If $(a_n)_{n\in\N}$ and $(b_n)_{n\in\N}$ are sequences in $D$ with $0\le a_n \le b_n$ for $n\in\N$ } \\
&\mbox{and $\sum_n b_n$ exists in $E$, then so does $\sum_n a_n$}.
\end{align}
It is not difficult to see that $D$ is splitting in $E$ if and only if 
\begin{align}
\label{eqn:splitting2}
\notag &\mbox{If $(a_n)_{n\in\N}$ is a sequence in $D^+$ and $\sum_n a_n$ exists in $E$,} \\
&\mbox{then so does $\sum_n \1_A(n) a_n$ for all $A\subset \N$}.
\end{align}
If $D$ is splitting in $E$, then so is every linear subspace of $D$. 
If $E$ is $\sigma$-Dedekind complete, then $E$ is also splitting. 
More generally, $D$ is splitting in $E$ if every bounded increasing sequence in $D$ has a supremum in $E$. 
Also, $\R^2$ with the lexicographical ordering is splitting. 
\end{definition}

In Theorem \ref{theorem:varphi_of_Gamma_R-closed_then_L_and_V_stable} we will see what is the use of this concept. First, we have a look at the connection between ``splitting'' and ``mediated''.

\begin{lemma}
\label{lemma:splitting_implies_sort_of_mediated}
Suppose $D$ is a linear subspace of $E$. 
Consider the condition:
\begin{align}
\label{eqn:sort_of_mediated}
\notag 
& \mbox{For all sequences $(a_n)_{n\in\N}, (b_n)_{n\in\N}$ in $D$: }\\
& a_n \downarrow, b_n \uparrow, \ \infn a_n- b_n=0 \ \ \Longrightarrow \ \ 
 \infn a_n = \supn b_n. 
\end{align}
(The infima and suprema in \eqref{eqn:sort_of_mediated} are to be taken in $E$.)
If $D$ is either splitting or mediated in $E$, then \eqref{eqn:sort_of_mediated} holds. 
Conversely, \eqref{eqn:sort_of_mediated} implies that $D$ is splitting if $D=E$,
whereas \eqref{eqn:sort_of_mediated} implies that $D$ is mediated in $E$ if $E$ is a Riesz space and $D$ is a Riesz subspace of $E$. 
\end{lemma}
\begin{proof}
It will be clear that mediatedness implies \eqref{eqn:sort_of_mediated} and vice versa if $E$ is a Riesz space and $D$ a Riesz subspace of $E$. \\
If $D$ is splitting in $E$ and $a_n\downarrow, b_n \uparrow$ and $\inf a_n - b_n =0$, then 
$\sum_n b_{n+1}- b_n +  a_{n} - a_{n+1} = a_1 - b_1$. Hence \eqref{eqn:sort_of_mediated} holds. \\
Suppose $D=E$ and \eqref{eqn:sort_of_mediated} holds. 
Let $(a_n)_{n\in\N}$ and $(b_n)_{n\in\N}$ be sequences in $D$ with $0\le a_n \le b_n$ for $n\in\N$ such that $\sum_n b_n$ exists. Let $z  = \sum_n b_n$, $A_n = \sum_{i=1}^n a_i$, $C_n = \sum_{i=1}^n b_i - a_i$ for $n\in\N$. Then $A_n \uparrow, C_n \uparrow$ and $z- C_n - A_n \downarrow 0$ (note that $z- C_n \in D$). Hence $\supn A_n = \sum_n a_n$ exists. 
\end{proof}

\begin{observation}
\label{observation:mediated_and_splitting}
\begin{enumerate}
[label=(\alph*),align=left,leftmargin=1cm,topsep=3pt,itemsep=0pt,itemindent=-1.2em,parsep=0pt]
\item If $E$ is a Riesz space, then every splitting Riesz subspace is mediated in $E$. 
\item If $E$ is mediated, then it is splitting. The converse is also true if $E$ is a Riesz space. 
\label{observation:part_mediated_implies_splitting_in_itself}
\item $c_{00}$ is mediated in $c$, not splitting in $c$ (with $D=E=c$ also \eqref{eqn:sort_of_mediated} is not satisfied)). 
\item If $D$ is the space of all polynomial functions on $[0,1]$ with degree at most $2$ and $E=C[0,1]$, then $D$ is splitting in $E$, but not mediated in $E$. (Actually, $D$ is splitting, but not mediated.)  \\
$D$ is splitting (and satisfies \eqref{eqn:sort_of_mediated} with $E=D$): If $u_n\in E^+$, $u_n \uparrow$ and $u_n \le \1$, then $|u_n(x)- u_n(y)|\le 4|x-y|$ as can be concluded from the postscript in Example \ref{example:gamma_LV_subsetneq_Gamma_VL}. Therefore the pointwise supremum is continuous. It is even in $D$ since $u_n(x) = a_n x^2 + b_n x + c_n$, where $a_n, b_n,c_n$ are linear combinations of $u_n(0),u_n(\frac12), u_n(1)$ (see also the postscript in Example \ref{example:gamma_LV_subsetneq_Gamma_VL}). \\
$D$ is not mediated: For example one can find countable $A,B \subset E$ for which $\1_{[\frac12,1])}$ is pointwise the infimum of $A$ and $\1_{(\frac12,1)}$ is pointwise the supremum of $B$, then $\inf A- B=0$, but there is no $h\in E$ with $B\le h \le A$.)
\end{enumerate}
\end{observation}

\begin{theorem}
\label{theorem:banach_lattice_sigma_order_c_norm_mediated_and_splitting}
Let $E$ be a Banach lattice with $\sigma$-order continuous norm. 
Then $E$ is both mediated and splitting. 
\end{theorem}
\begin{proof}
Suppose $a_n,b_n\in E$ with $0 \le a_n \le b_n$ for $n\in\N$. 
Suppose that $\{\sum_{n=1}^N b_n: N\in\N\}$ has a supremum $s$ in $E$. 
We prove that $\{\sum_{n=1}^N a_n: N\in\N\}$ has a supremum in $E$. 
Since the norm is $\sigma$-order continuous, we have $\|s- \sum_{n=1}^N b_n\| \rightarrow 0$. 
In particular we get that for all $\epsilon>0$ there exists an $N\in\N$ such that for all $n,m\ge N$ with $m>n$ we have $\|\sum_{i=n}^m b_i\|<\epsilon$ and thus $\|\sum_{i=n}^m a_i\|<\epsilon$. From this we infer that $(\sum_{n=1}^N a_n)_{N\in\N}$ converges in norm. 
Therefore it has a supremum in $E$. Thus $E$ is splitting. 
By Lemma \ref{lemma:splitting_implies_sort_of_mediated} $E$ is mediated. 
\end{proof}

\begin{theorem} \mbox{}
\label{theorem:varphi_of_Gamma_R-closed_then_L_and_V_stable}
\begin{enumerate}[label=\emph{(\alph*)},align=left,leftmargin=1cm,topsep=3pt,itemsep=0pt,itemindent=-1.2em,parsep=0pt]
\item 
$\varphi(\Gamma)$ splitting in $E$
$\Longrightarrow $ $\Gamma_L$ is stable and $\varphi_L$ is laterally extendable. 
\item $\varphi(\Gamma)$ mediated in $E$
$\Longrightarrow $ $\Gamma_V$ is stable and $\varphi_V$ is laterally extendable. 
\item $\varphi(\Gamma)$ splitting in $E$ and $\varphi_L(\Gamma_L)$ mediated in $E$
$\Longrightarrow $ $(\Gamma_{L})_V$ is stable and $(\varphi_{L})_V$ is laterally extendable. 
\end{enumerate}
\end{theorem}
\begin{proof}
(a) Let $f\in \Gamma_L$, $B\in \cI$; we prove $f\1_B \in \Gamma_L$. (This is sufficient by Theorem \ref{theorem:stability_implies_laterally_extendability_for_V_and_L}(a).) 
Without loss of generality, assume $f\ge 0$. 
Choose a $\varphi$-partition  $(A_n)_{n\in\N}$  for $f$. 
Now apply \eqref{eqn:splitting} to
\begin{align}
 a_n := \varphi( f \1_{A_n \cap B}), \quad b_n:= \varphi(f\1_{A_n}) \qquad (n\in\N). 
\end{align} 
(b) follows from Lemma \ref{lemma:bounded_measurable_functions} and Theorem \ref{theorem:stability_implies_laterally_extendability_for_V_and_L}(b). \\
(c) By (a) $\Gamma_L$ is stable and $\varphi_L$ is laterally extendable. Hence we can apply (b) to $\Gamma_L$ and $\varphi_L$ (instead of $\Gamma$ and $\varphi$) and obtain (c). 
\end{proof}

\begin{observation}
\label{observation:in_canonical_example_E_is_complete_iff_Gamma_L_is_stable}
To some extent, the assumption 
of Theorem \ref{theorem:varphi_of_Gamma_R-closed_then_L_and_V_stable}(a) is minimal. \\
Indeed, in the situation of Example \ref{example:example_latterally_extendable_sequences}, we see that $\Gamma_L$ is stable if and only if $E$ (which is $\varphi(\Gamma)$) is splitting in $F$ (see \eqref{eqn:splitting2}).
\end{observation}

In Theorem \ref{theorem:varphi_of_Gamma_R-closed_then_L_and_V_stable}(c) we assumed that $\varphi_L(\Gamma_L)$ (and thus also $\varphi(\Gamma)$) was mediated in $E$. 
It may happen that $\varphi(\Gamma)$ is mediated in $E$, but $\varphi_L(\Gamma_L)$ is not, as Example \ref{example:varphi_of_Gamma_L_not_mediated} illustrates. 
However, splitting is preserved under the lateral extension and mediation is preserved under the vertical extension, see Theorem \ref{theorem:closability_preserved}. 

\begin{example}
\label{example:varphi_of_Gamma_L_not_mediated}
Let $X=\N$, $\cI= \cP(\N)$, $E=F=c$. 
Let $\Gamma = c_{00}[c_{00}] $ (see \S\ref{section:some_notation}) and $\varphi: \Gamma \rightarrow E$ be given by $\varphi(f) = \sum_{n\in\N} f(n)$. 
Then $\varphi(\Gamma)= c_{00}$, which is mediated in $c$.
A function $f: \N \rightarrow c$ is partially in $\Gamma$ if and only if $f(\N) \subset c_{00}$. 
For $x\in c^+$ the function given by $f(n) = x(n) \1_{\{n\}}$ for $n\in\N$ lies in $\Gamma_L$, and $\varphi_L(f) =x$. It follows that $\varphi_L(\Gamma_L)$ is $c$, which is not mediated in $c$. 
\end{example}

\begin{theorem} \
\label{theorem:closability_preserved}
\begin{enumerate}[label=\emph{(\alph*)},align=left,leftmargin=1cm,topsep=3pt,itemsep=0pt,itemindent=-1.2em,parsep=0pt]
\item
If $\varphi(\Gamma)$ is splitting in $E$, then so is  $\varphi_L(\Gamma_L)$. 
\item 
If $\varphi(\Gamma)$ is mediated in $E$, then so is  $\varphi_V(\Gamma_V)$. 
\end{enumerate} 
\end{theorem}
\begin{proof}
(a) Suppose $a_n\in \varphi_L(\Gamma_L)^+$ for $n\in\N$ and $\sum_n a_n$ exists. 
Let $A\subset \N$. 
For all $n\in\N$ there exist $b_{n1},b_{n2},\dots \in \varphi(\Gamma)^+$ with $a_n = \sum_m b_{nm}$. 
Hence $\sum_n a_n = \sum_{n,m} b_{nm}$ and so $\sum_{n,m} \1_{A\times \N}(n,m) b_{nm} = \sum_n \1_A(n) a_n$ exists in $E$. 

(b) Suppose $A,B\subset \varphi_V(\Gamma_V)$ are countable sets with $\inf A-B=0$. For all $a\in A$ and $b\in B$ there exist countable sets $\Upsilon_a,\Lambda_b\subset \Gamma$ with $a= \inf \varphi(\Upsilon_a), b= \sup \varphi(\Lambda_b)$. 
Then $\inf \varphi ( \bigcup_{a\in A} \Upsilon_a - \bigcup_{b\in B} \Lambda_b)=0$ and thus $\inf A = \inf \varphi  (\bigcup_{a\in A} \Upsilon_a )= \sup \varphi( \bigcup_{b\in B} \Lambda_b)= \sup B$. 
\end{proof}

\begin{observation}
For a Riesz space $F$ we will now investigate under which conditions the space $\Gamma_L$ is a Riesz subspace of $F^X$.
The next example shows that even if $E$ is a Riesz space and $\Gamma$ is a Riesz subspace of $F^X$, $\Gamma_L$ may not be one. However, see Theorem \ref{theorem:splitting_E_implies_Gamma_L_Riesz}.
\end{observation}

\begin{example}
\label{example:lateral_extension_not_Riesz}
Let $a,b$ be as in Example \ref{example:gamma_L_not_stable}; this time put $d=(0,a_1+a_2,0,a_3+a_4,\dots)$. 
Then $a,d\in \Gamma_L$ but $a\wedge d = b \notin \Gamma_L$. 

Hence, in Example 
\ref{example:example_latterally_extendable_sequences},
if $F$ is a Riesz space and $E$ is not splitting in $F$, 
then $\Gamma_L$ is not a Riesz subspace of $F^X$. 
As we will see in Theorem \ref{theorem:splitting_E_implies_Gamma_L_Riesz}, considering the situation of Example \ref{example:example_latterally_extendable_sequences}: 
$\Gamma_L$ is a Riesz subspace of $F^X$ if and only if $E$ is splitting in $F$. 
\end{example}

\begin{lemma} 
\label{lemma:splitting_implies_Gamma_L_Riesz_ideal}
Let $f: X \rightarrow F$ be partially in $\Gamma$. 
\begin{enumerate}[label=\emph{(\alph*)},align=left,leftmargin=1cm,topsep=3pt,itemsep=0pt,itemindent=-1.2em,parsep=0pt]
\item 
If $f$ is in $ \Gamma_{LV}$, then $f \in \Gamma_L$. 
\item Suppose $\varphi(\Gamma)$ is splitting in $E$. 
If $g\le f \le h$ for certain $g,h\in \Gamma_L$, then $f\in \Gamma_L$. 
\end{enumerate}
\end{lemma}
\begin{proof}
(a) By the definition of $\Gamma_{LV}$ there exists a $\rho \in \Gamma_L$ with $\rho \le f$.
Then $ f- \rho$ is partially in $\Gamma$, $f-\rho \in \Gamma_{LV}$, and we are done if $f- \rho \in \Gamma_L$. 
Hence we may assume $f \ge 0$. \\
Let $(A_n)_{n\in\N}$ be a partition for $f$; we prove $\sum_n \varphi(f\1_{A_n}) = \varphi_{LV}(f)$. 
It will be clear that $\sum_{n=1}^N \varphi(f \1_{A_n}) \le  \varphi_{LV} ( f) $ for $N\in\N$. 
For the reverse inequality let $h\in E$ be an upper bound for $\{ \sum_{n=1}^N \varphi(f \1_{A_n}): N\in\N\}$. 
It suffices to show that $h$ must be an upper bound for $\{ \varphi_L(\sigma): \sigma \in \Gamma_L, \sigma \le f\}$. \\
Take a $\sigma \in \Gamma_L$ with $\sigma \le f$. 
If $(B_n)_{n\in\N}$ is any refinement of $(A_n)_{n\in\N}$ that is a $\varphi$-partition for $\sigma$, then for all $M\in\N$ there exists an $N\in\N$ with $B_1 \cup \cdots \cup B_M \subset A_1 \cup \cdots \cup A_N$, so that $h \ge \sum_{n=1}^N \varphi(f \1_{A_n}) \ge  \sum_{m=1}^M \varphi(f \1_{B_m})\ge  \sum_{m=1}^M \varphi(\sigma \1_{B_m})$.
It follows from Lemma \ref{lemma:a_tool_for_refinements}, applied to $\sigma$, that the partition $(B_m)_{m\in\N}$ can be chosen so that this implies $h \ge \varphi_L(\sigma)$.  

(b) As $h-g \in \Gamma_L$ and $0 \le f-g \le h-g$, we may (and do) assume $g=0$. 
Let $(A_n)_{n\in\N}$ be a partition for $f$ that is also a $\varphi$-partition for $h$. 
Now just apply \eqref{eqn:splitting} to
\begin{align}
 a_n:= \varphi(f \1_{A_n}), \quad b_n:=  \varphi(h\1_{A_n}) \qquad (n\in\N). 
\end{align}
\
\end{proof}

As a consequence of Lemma \ref{lemma:splitting_implies_Gamma_L_Riesz_ideal}:

\begin{theorem}
\label{theorem:splitting_E_implies_Gamma_L_Riesz}
Let $F$ be a Riesz space and $\Gamma$ be a Riesz subspace of $F^X$. 
The functions $X\rightarrow F$ that are partially in $\Gamma$ form a Riesz space, $\Xi$.
If $\varphi(\Gamma)$ is splitting in $E$, then $\Gamma_L$ is a Riesz ideal in $\Xi$, in particular, $\Gamma_L$ is a Riesz space. 
\end{theorem}

In the classical integration theory and the Bochner integration theory 
one starts with considering a measure space $(X,\cA,\mu)$ and simple functions on $X$ with values in $\R$ or in a Banach space. 
One defines an integral on these simple functions using the measure and extends this integral to a larger class of integrable functions. 
In \ref{obs:simple_functions_and_their_integral} we will 
follow a similar procedure, replacing $\R$ or the Banach space with $E$ and applying the lateral extension. 
In Section \ref{section:extensions_of_simple_functions} we will treat such extensions in more detail.

\begin{observation}
\label{obs:simple_functions_and_their_integral}
Suppose $(X,\cA,\mu)$ is a $\sigma$-finite complete measure space 
and suppose $E$ is directed. Let $F=E$. 
For $\cI$ we choose $\{A\in \cA: \mu(A)<\infty\}$.
The $\sigma$-finiteness of $\mu$ guarantees the existence of a partition (and vice versa). 

We say that a function $f: X\rightarrow E$ is \emph{simple} if there exist $N\in\N$, $a_1,\dots,a_N\in E$, $A_1,\dots, A_N\in \cI$ for which 
\begin{align}
\label{eqn:representation_simple_function}
f = \sum_{n=1}^N a_n \1_{A_n}.
\end{align}
The simple functions form a stable directed linear subspace $S$ of $E^X$, which is a Riesz subspace of $E^X$ in case $E$ is a Riesz space. 

For a given $f$ in $S$ one can choose a representation \eqref{eqn:representation_simple_function} in which the sets $A_1,\dots, A_N$ are pairwise disjoint; thanks to the $\sigma$-finiteness of $\mu$ one can choose them in such a way that they occur in a partition $(A_n)_{n\in\N}$.

This $S$ is going to be our $\Gamma$. We define $\varphi : S\rightarrow E$ by 
\begin{align}
\label{eqn:integral_for_simple_function}
\varphi(f) = \sum_{n=1}^N \mu(A_n) a_n,
\end{align}
where $f,N,A_n,a_n$ are as in \eqref{eqn:representation_simple_function}. 
The $\sigma$-additivity of $\mu$ is (necessary and) sufficient to show that $S$ is laterally extendable. 

A function $f: X\rightarrow E$ is partially in $S$ if and only if there exist a partition $(A_n)_{n\in\N}$ and a sequence $(a_n)_{n\in\N}$ in $E$ for which 
\begin{align}
\label{eqn:partially_simple}
f = \sum_{n\in\N} a_n \1_{A_n}. 
\end{align}
An $f$ as in \eqref{eqn:partially_simple} with $f\ge 0$  that is partially in $S$ is an element of $S_L$ if and only if $\sum_n \mu(A_n) a_n$ exists in $E$. (See Theorem \ref{theorem:unique_integral_elements_of_gamma_sigma}.)
\end{observation}

\section{Combining vertical and lateral extensions} 
\label{section:combinations_of_extensions}

\textbf{\emph{In this section $E,F,X,\cI,\Gamma,\varphi$ are as in Section \ref{section:lateral_extension}.}} 

As we have seen, the lateral extension differs from the vertical extension in the sense that the vertical extensions of $\Gamma$ and $\varphi$ can always be made, but for lateral extension we had to assume the space $\Gamma$  to be stable and $\varphi$ to be laterally extendable (see \ref{observation:lateral_integral_well_defined}). 
In this section we investigate when one can make a lateral extension of another (say vertical) extension. Furthermore we will compare different extensions and combinations of extensions. 

\bigskip

Instead of $(\Gamma_L)_V$ and $((\Gamma_L)_V)_L$ we write $\Gamma_{LV}$ and $\Gamma_{LVL}$; similarly $\varphi_{LV}= (\varphi_{L})_V$ etc.

\begin{observation}
By Theorem \ref{theorem:stability_implies_laterally_extendability_for_V_and_L}  the following holds for a stable directed linear subspace $\Delta$ of $F^X$ and a laterally extendable  order preserving linear map $\omega: \Delta \rightarrow E$: 
If $\Delta_L$ is stable, then $\omega_L$ is laterally extendable (and so $\Delta_{LL}$ exists). 
If $\Delta_V$ is stable, then $\omega_V$ is laterally extendable (and so $\Delta_{VL}$ exists). We will use these facts without explicit mention. 
\end{observation}

\begin{observation}
\label{observation:integrals_carry_over_to_larger_spaces}
The following statements follow from the definitions and theorems we have:
 \newcounter{integrals_carry_over}
\begin{enumerate}[label=(\alph*),align=left,leftmargin=1cm,topsep=3pt,itemsep=0pt,itemindent=-1.2em,parsep=0pt]
\item 
$\Gamma_V \subset \Gamma_{LV}$ and $\varphi_{LV}=\varphi_V$ on $\Gamma_V$. 
\item 
$\Gamma_L \subset \Gamma_{LV}$ and $\varphi_{LV}=\varphi_L$ on $\Gamma_L$. 
\item $\varphi_V= \varphi_L$ on $\Gamma_L\cap \Gamma_V$.
\setcounter{integrals_carry_over}{\value{enumi}}
\end{enumerate}
For (d), (e) and (f) let $\Gamma_V$ be stable.
\begin{enumerate}[label=(\alph*),align=left,leftmargin=1cm,topsep=3pt,itemsep=0pt,itemindent=-1.2em,parsep=0pt]
\setcounter{enumi}{\value{integrals_carry_over}}
\item 
$\Gamma_{LV} \subset \Gamma_{VLV}$ and $\varphi_{VLV}=\varphi_{LV}$ on $\Gamma_{LV}$. 
\item 
$\Gamma_{VL} \subset \Gamma_{VLV}$ and $\varphi_{VLV} = \varphi_{VL}$ on $\Gamma_{VL}$. 
\item 
$\varphi_{LV}= \varphi_{VL}$ on $\Gamma_{LV} \cap \Gamma_{VL}$. 
\end{enumerate}
Observe that as a consequence of (a) and (b): 
If $f\in \Gamma_L$ and $g\in \Gamma_V$ and $f\le g$ (or $f\ge g$), then $\varphi_L(f) \le \varphi_V(g)$ (or $\varphi_L(f) \ge \varphi_V(g)$). 
Moreover, as a consequence of (c) and (d); if $\Gamma_V$ is stable: 
If $f\in \Gamma_{LV}$ and $g\in \Gamma_{VLV}$ and $f\le g$ (or $f\ge g$), then $\varphi_{LV}(f) \le \varphi_{VLV}(g)$ (or $\varphi_{LV}(f) \ge \varphi_{VLV}(g)$). 
\end{observation}

\begin{observation}
\label{observation:some_remarks_on_biggest_extension}
Note that if $\Gamma$ is stable and $\varphi$ is laterally extendable, then we can extend $\Gamma$ to $\Gamma_V, \Gamma_L$ and $\Gamma_{LV}$. 
If, moreover, $\Gamma_V$ is stable, then we can also extend $\Gamma$ to $\Gamma_{VL}$ and $\Gamma_{VLV}$. 
However, ``more stability'' will not give us larger extensions than $\Gamma_{VLV}$. Indeed, 
if $\Gamma_{LV}$ is stable then $\Gamma_{LV} \subset \Gamma_{LVL}=\Gamma_{VL}$ (see Theorem \ref{theorem:Gamma_LV_is_subset_of_Gamma_VL}). 
If moreover $\Gamma_{VLV}$ is stable, then even $\Gamma_{VLVL}= \Gamma_{VLV}=\Gamma_{VL}$. 
\end{observation}

\begin{lemma} \mbox{}
\label{lemma:positive_integrable_functions_and_countable_minorizers_in_Gamma}
\begin{enumerate}[label=\emph{(\alph*)},align=left,leftmargin=1cm,topsep=3pt,itemsep=0pt,itemindent=-1.2em,parsep=0pt]
\item If $f\in \Gamma_{LV}^+$, then there exists a countable $\Lambda \subset \Gamma$ with $\Lambda \le f$ and $ \varphi_{LV}(f) = \sup \varphi(\Lambda)$. 
\item If $\Gamma_V$ is stable and $f\in \Gamma_{VL}^+$, then there exists a countable $\Lambda  \subset \Gamma$ with $\Lambda \le f$ and $ \varphi_{VL}(f) = \sup \varphi(\Lambda)$. 
\end{enumerate}
\end{lemma}
\begin{proof}
(a) There exist $\sigma_1,\sigma_2, \dots$ in $\Gamma_L$ with $\sigma_n \le f$ for all $n\in\N$ and $\supn \varphi_{L}(\sigma_n) = \varphi_{LV}(f)$. Hence, we are done if for every $\sigma$ in $\Gamma_L$ with $\sigma \le f$ there is a countable set $\Lambda_\sigma \subset \{ \rho\in \Gamma: \rho \le f\}$ such that every upper bound for $\varphi(\Lambda_\sigma)$ majorizes $\varphi_L(\sigma)$. But that is not hard to prove. 
For such a $\sigma$, by Lemma \ref{lemma:a_tool_for_refinements} there exists a partition $(B_m)_{m\in\N}$ for which \eqref{eqn:tool_refinements_property} holds. 
Now let $\Lambda_\sigma$ be $\{\sum_{m=1}^M \sigma \1_{B_m}: M\in\N\}$. \\
(b) Suppose  $\Gamma_V$ is stable. 
Let $(A_n)_{n\in\N}$ be a $\varphi_V$-partition for $f$. 
Then the set $\Lambda_f=\{ \sum_{n=1}^N f \1_{A_n}: N\in\N\}$ is a countable subset of $\Gamma_V$ and $\sup \varphi_V(\Lambda_f) = \varphi_{VL}(f)$. Moreover, for every $N\in\N$ there is a countable set $\Lambda_N \subset \{ \sigma \in \Gamma: \sigma \le \sum_{n=1}^N f\1_{A_n}\}$ for which $\sup \varphi(\Lambda_N) = \varphi_{V}( \sum_{n=1}^N f\1_{A_n})$. 
Take $\Lambda = \bigcup_{N\in\N} \Lambda_N $. 
\end{proof}

\begin{theorem}
\label{theorem:extended_functions_inbetween_Gamma_functions}
For (b),(c),(d) and (e)  let $\Gamma_V$ be stable and $f$ be partially in $\Gamma_V$. 
\begin{enumerate}[label=\emph{(\alph*)},align=left,leftmargin=1cm,topsep=3pt,itemsep=0pt,itemindent=-1.2em,parsep=0pt]
\item 
\label{item:if_f_in_LV_then_f_in_V_iff_inclosed_by_Gamma}
If $f\in \Gamma_{LV}$, then
\begin{align*}
f\in \Gamma_V\iff \mbox{ there exist }\pi,\rho \in \Gamma \mbox{ with }\pi \le f \le \rho \  \footnotemark. 
\end{align*}
\item
 If $f\in \Gamma_{VL}$, then
\begin{align*}
f\in \Gamma_V\iff \mbox{ there exist }\pi,\rho \in \Gamma \mbox{ with }\pi \le f \le \rho \ 
\textsuperscript{\emph{\ref{footnote:gamma_V_ideal_in_LV}}}.
\end{align*}
\item 
\label{item:f_in_LV_iff_in_VL_and_inclosed_by_L}
\quad $
f \in \Gamma_{LV}  \iff 
f\in \Gamma_{VL} \mbox{ and there exist } \pi,\rho \in \Gamma_L \mbox{ with } \pi \le f \le \rho. 
$ \vspace{2mm}
\item If $\varphi_V(\Gamma_V)$ is splitting in $E$, then 
\label{item:splitting_then_f_in_VL_iff_f_inclosed_by_VL}
\begin{align*}
f\in \Gamma_{VL} \iff \mbox{ there exist } \pi,\rho \in \Gamma_{VL} \mbox{ with } \pi \le f \le \rho. 
\end{align*}
\item 
\label{item:splitting_then_f_in_VL_and_LV_iff_f_inclosed_by_L}
If $\varphi_V(\Gamma_V)$ is splitting in $E$, then
\begin{align*}
f\in \Gamma_{VL}\cap \Gamma_{LV} \iff \mbox{ there exist } \pi,\rho \in \Gamma_{L} \mbox{ with } \pi \le f \le \rho. 
\end{align*}
\end{enumerate}
\footnotetext{\label{footnote:gamma_V_ideal_in_LV}By the definition of ideal in \cite{Bo54} or  \cite{fuchs1966riesz} (note that $\Gamma_V$ is directed)
 this means that $\Gamma_V$ is the smallest ideal in $\Gamma_{LV}$ (and for (b); in $\Gamma_{VL}$) that contains $\Gamma$.}
\end{theorem}
\begin{proof}
The proofs of (a) and (b) are similar to the proof of (c) and therefore omitted. 

(c) $\Leftarrow$:
By Lemma \ref{lemma:positive_integrable_functions_and_countable_minorizers_in_Gamma} (b) there exist countable sets $\Lambda, \Upsilon \subset \Gamma$ with $\Lambda \le f- \pi$ and $\Upsilon \le \rho - f $  for which $\sup \varphi(\Lambda) = \varphi_{VL}(f-\pi)$ and $\sup \varphi(\Upsilon) = \varphi_{VL}(\rho - f)$. 
Then $\Lambda + \pi$ and $\rho - \Upsilon$ are countable subsets of $\Gamma_L$ with $\Lambda + \pi \le f \le \rho - \Upsilon$ and $\sup \varphi_L(\Lambda + \pi) = \varphi_{VL}(f) = \inf \varphi_L (\rho - \Upsilon)$. 
Hence $f\in \Gamma_{LV}$. 

$\Rightarrow$: Let $f\in \Gamma_{LV}$ and be partially in $\Gamma_V$. 
There exists a $\pi \in \Gamma_L$ for which $f- \pi \in \Gamma_{LV}^+$, hence we may assume $f\ge 0$. 
Let $(A_n)_{n\in\N}$ be a $\Gamma_V$-partition for $f$, i.e., 
$f\1_{A_n}\in \Gamma_V$ and thus $\varphi_{LV}(f\1_{A_n}) = \varphi_V(f\1_{A_n})$ for all $n\in\N$
(see \ref{observation:integrals_carry_over_to_larger_spaces}(a)). 
Then $\varphi_{LV}(f) \ge \sum_{n=1}^N \varphi_V(f\1_{A_n})$ for all $N\in\N$.
Let $h\in E$ be such that $h \ge \sum_{n=1}^N \varphi_V(f\1_{A_n})$ for all $N\in\N$. 
From Lemma \ref{lemma:a_tool_for_refinements} we infer that $h\ge \varphi_L(\sigma)$ for every $\sigma \in \Gamma_L$ with $\sigma \le f$. 
We conclude that $\sum_n \varphi_V(f\1_{A_n}) = \varphi_{LV}(f)$, i.e., $f\in \Gamma_{VL}$.

(d) $\Leftarrow$: We may assume $\pi =0$. 
Let $(A_n)_{n\in\N}$ be a $\varphi_V$-partition for $\rho$ with  $f\1_{A_n}\in \Gamma_V$ for all $n\in\N$. Then $0\le \varphi_V(f\1_{A_n}) \le \varphi_V(\rho \1_{A_n})$ for all $n\in\N$ and $\sum_{n} \varphi_V(\rho\1_{A_n})$ exists in $E$. Hence, so does $\sum_{n} \varphi_V(f \1_{A_n})$, i.e., $f\in \Gamma_{VL}$. 

(e) is a consequence of (c) and (d). 
\end{proof}

In the following example all functions in $\Gamma_{LV}$ are partially in $\Gamma_V$. 

\begin{example}
\label{example:elements_of_gammaLV_are_partially_in_gammaV}
Consider $X=\N, \cI=\cP(\N)$, $E=F$; let $D$ be a linear subspace of  $E$ and let $D_V$ be the vertical extension of $D$ with respect to the inclusion map $D\rightarrow E$. 
Let $\Gamma = c_{00}[D]$ and $\varphi: \Gamma \rightarrow E$ be $\varphi(f) = \sum_{n\in\N}f(n)$. 
Then $\Gamma_V = c_{00}[D_V]$. 
Let $f\in \Gamma_{LV}$. We will show that $f(k) \in D_V$ and thus that $f$ is partially in $\Gamma_V$. 
Let $\sigma_n, \tau_n \in \Gamma_L$ be such that $\sigma_n \le  f \le  \tau_n$ and $\infn \varphi(\tau_n) = \supn \varphi(\sigma_n)$. Then $\infn (\tau_n(k) - \sigma_n(k)) \le \infn \varphi(\tau_n - \sigma_n) = 0$. Since $ \sigma_n(k), \tau_n(k) \in D$ for all $n\in\N$, we have $f(k) \in D_V$. \\
Thus every $f\in \Gamma_{LV}$ is partially in $\Gamma_V$. 
Since $\Gamma_V$ is stable, by Theorem \ref{theorem:extended_functions_inbetween_Gamma_functions}\emph{\ref{item:f_in_LV_iff_in_VL_and_inclosed_by_L}} we conclude that $\Gamma_{LV} \subset \Gamma_{VL}$.
\end{example}

\begin{lemma}
\label{lemma:gamma_LV_stable_then_every_function_of_it_is_partially_V}
Suppose that $\Gamma_{LV}$ is stable. 
Then every $f\in \Gamma_{LV}$ is partially in $\Gamma_V$. 
\end{lemma}
\begin{proof}
Let $f\in \Gamma_{LV}$ and let $\pi,\rho\in \Gamma_L$ be such that $\pi \le f \le \rho$. 
Let $(A_n)_{n\in\N}$ be a $\varphi$-partition for both $\pi$ and $\rho$. 
Then $f\1_{A_n} \in \Gamma_{LV}$ and $\pi \1_{A_n} \le f \1_{A_n} \le \rho \1_{A_n}$ for all $n\in\N$.
By Theorem \ref{theorem:extended_functions_inbetween_Gamma_functions}\emph{\ref{item:if_f_in_LV_then_f_in_V_iff_inclosed_by_Gamma}} we conclude that $f\1_{A_n}\in \Gamma_V$.
\end{proof}

\begin{theorem}
\label{theorem:Gamma_LV_is_subset_of_Gamma_VL}
Suppose that $\Gamma_V$ and $\Gamma_{LV}$ are stable. 
Then $\Gamma_{LV} \subset \Gamma_{VL}=\Gamma_{LVL}$. 
Write $\overline \Gamma = \Gamma_{VL}$ and $\overline \varphi = \varphi_{VL}$. 
If $\overline \Gamma$ is stable, then $\overline \Gamma_L = \overline \Gamma$ and $\overline \varphi_L = \overline \varphi$.
If $\overline \Gamma_{V}$ is stable, then $\overline \Gamma_V = \overline \Gamma$ and 
$\overline \varphi_V = \overline \varphi$.  \\
In particular, if $\varphi_L(\Gamma_L)$ is mediated in $E$ and $\varphi_V(\Gamma_V)$ is splitting in $E$, then $\Gamma_V$, $\Gamma_{LV}$ and $\Gamma_{VL}$ are stable (see Theorem \ref{theorem:varphi_of_Gamma_R-closed_then_L_and_V_stable}) and thus $\Gamma_{LV}\subset \overline \Gamma$, $\overline \Gamma = \overline \Gamma_V = \overline \Gamma_L$, $\overline \varphi = \overline \varphi_V = \overline \varphi_L$, so $\overline{\overline \Gamma} = \overline \Gamma$ (and $\overline{ \overline \varphi } = \overline \varphi$). 
\end{theorem}
\begin{proof}
The inclusion $\Gamma_{LV}\subset \Gamma_{VL}$ follows by Theorem \ref{theorem:extended_functions_inbetween_Gamma_functions}\emph{\ref{item:f_in_LV_iff_in_VL_and_inclosed_by_L}} and Lemma \ref{lemma:gamma_LV_stable_then_every_function_of_it_is_partially_V}. We prove $\Gamma_{LVL}\subset \Gamma_{VL}$.
For $f\in \Gamma_{LVL}^+$ there is a $\varphi_{LV}$-partition for $f$ and since $\Gamma_{LV} \subset \Gamma_{VL}$ this is also a $\varphi_{VL}$-partition for $f$, hence there exists a $\varphi_V$-partition for $f$, i.e., $f\in \Gamma_{VL}$. 

Suppose $\overline \Gamma$ is stable. Then $\overline \Gamma_{L} = (\Gamma_{VL})_L = \Gamma_{VL}= \overline \Gamma$ and $\overline \varphi_L = \overline \varphi$ by Theorem \ref{theorem:stability_implies_laterally_extendability_for_V_and_L}(a). 

Suppose $\overline \Gamma_V$ to be stable. 
As $\Gamma_V$ is stable we can apply the first part of the theorem to $\Gamma_V$ instead of $\Gamma$. 
Indeed, $(\Gamma_V)_V$ and $(\Gamma_V)_{LV}$ are stable, since $(\Gamma_V)_V=\Gamma_V$ and $(\Gamma_V)_{LV}= \overline \Gamma_V$. 
Hence, $(\Gamma_V)_{LV} \subset (\Gamma_V)_{VL} = \Gamma_{VL}$, i.e., $\overline \Gamma_V \subset \overline \Gamma$ (and $\overline \varphi_V = \overline \varphi$).

Suppose $\varphi_L(\Gamma_L)$ is mediated in $E$ and $\varphi_V(\Gamma_V)$ is splitting in $E$.
Then $\Gamma_L$, $\Gamma_V$ and $\Gamma_{LV}$ are stable by Theorem \ref{theorem:varphi_of_Gamma_R-closed_then_L_and_V_stable}(a),(b) and (c). 
Consequently, again by Theorem \ref{theorem:varphi_of_Gamma_R-closed_then_L_and_V_stable}(b) $\Gamma_{VL}$ is stable.
\end{proof}

\begin{corollary}
Suppose $E$ is mediated (and thus splitting), $\overline \Gamma = \Gamma_{VL}$. 
Then $\overline \Gamma = \overline \Gamma_V = \overline \Gamma_L$, so $\overline{\overline \Gamma} = \overline \Gamma$ (and $\overline{ \overline \varphi } = \overline \varphi$). 
\end{corollary}

At the end of \S\ref{section:combinations_of_extensions} we will 
show that sometimes $\Gamma_{VL} \subsetneq \Gamma_{LV}$ (Example \ref{example:gamma_vl_subsetneq_gamma_lv})
and sometimes $\Gamma_{LV}\subsetneq \Gamma_{VL}$ (Example \ref{example:gamma_LV_subsetneq_Gamma_VL}). 
Note that this implies that $\Gamma_{VLV}$ can be strictly larger then either $\Gamma_{VL}$ or $\Gamma_{LV}$.

Theorem \ref{theorem:Gamma_LV_is_subset_of_Gamma_VL} raises the question whether stability of $\Gamma_V$ entails $\Gamma_{VL} \subset \Gamma_{LV}$. 
In general the answer is negative; see Example \ref{example:gamma_LV_subsetneq_Gamma_VL}.
In Theorem \ref{theorem:conditions_for_VL_subset_LV} we give conditions sufficient for the inclusion. 

\begin{theorem}
\label{theorem:conditions_for_VL_subset_LV}
Suppose $\Gamma_V$ is stable. Consider these two statements.
\begin{enumerate}[label=\emph{(\alph*)},align=left,leftmargin=1cm,topsep=3pt,itemsep=0pt,itemindent=-1.2em,parsep=0pt]
\item For every $f\in \Gamma_{VL}^+$ there is a $\rho$ in $\Gamma_L^+$ with $f\le \rho$. 
\item $E$ satisfies:
\begin{align}
\notag  &\mbox{If } Y_1,Y_2,\dots \subset E \mbox{ are nonempty countable with }\inf Y_n=0 \mbox{ for all } n\in\N, \\ 
&\mbox{then there exist } y_1\in Y_1,y_2 \in Y_2, \dots \mbox{ such that } \sum_n y_n \mbox{ exists in } E. 
\label{eqn:countable_sum_property}
\end{align}
\end{enumerate}
If (a) is satisfied, then $\Gamma_{VL}\subset \Gamma_{LV}$. 
(b) implies (a). 
\end{theorem}
\begin{proof}
If (a) is satisfied, then by Theorem \ref{theorem:extended_functions_inbetween_Gamma_functions}\emph{\ref{item:f_in_LV_iff_in_VL_and_inclosed_by_L}} follows that $\Gamma_{VL} \subset \Gamma_{LV}$. \\
Suppose (b). 
Let $f\in \Gamma_{VL}^+$.
Let $(A_n)_{n\in\N}$ be a $\varphi_V$-partition for $f$. 
For $n\in\N$, let $\Upsilon_n \subset \Gamma$ be a countable set with $f\1_{A_n} \le \Upsilon_n$ and 
\begin{align}
 \varphi_V(f\1_{A_n}) = \inf \varphi(\Upsilon_n). 
\end{align}
We may assume $\sigma \1_{A_n}= \sigma $ for all $\sigma \in \Upsilon_n$. 
Choose $\sigma_n \in \Upsilon_n$ for $n\in\N$ such that $\sum_n (\varphi(\sigma_n) - \varphi_V(f\1_{A_n}))$ and thus $\sum_n \varphi(\sigma_n)$ exist in $E$. 
Then $\rho :=  \sum_{n\in\N} \sigma_n$ is in $\Gamma_L^+$ with $f \le \rho $. 
\end{proof}

\begin{observation}
\label{observation:countable_sum_property}
We will discuss examples of spaces $E$ for which \eqref{eqn:countable_sum_property} holds. \\
(I) If $E$ is a Banach lattice with $\sigma$-order continuous norm, then $E$ satisfies \eqref{eqn:countable_sum_property} (one can find $y_n \in Y_n$ with $\|y_n\|\le 2^{-n}$). \\
(II)
Let $(X,\cA,\mu)$ be a  complete $\sigma$-finite  measure space and assume there exists a $g\in L^1(\mu)$ with $g>0$ $\mu$-a.e.. 
Then the space $E$ of equivalence classes of measurable functions $X\rightarrow \R$ satisfies \eqref{eqn:countable_sum_property}:
It is sufficient to prove that if $Z_1,Z_2,\dots \subset E$ are nonempty countable with $\inf Z_n=0$ for all $n\in\N$, then there exists $z_1\in Z_1,z_2\in Z_2,\dots$ and a $z\in E$ such that $z_n \le z$ for all $n\in\N$ (for $Z_n$ take $2^n Y_n$). 
One can prove that such a $z$ exists by mapping the equivalence classes of measurable functions into $L^1(\mu)$ by the order isomorphism $f\mapsto (\arctan \circ f)g$. \\
(III) $\R^\N$ is a special case of (II), therefore satisfies \eqref{eqn:countable_sum_property}.
\end{observation}

\begin{theorem}
\label{theorem:Gamma_VL=LV_for_banach_lattice_with_order_cont_norm}
Let $E$ be mediated and splitting and satisfy \eqref{eqn:countable_sum_property} (e.g. $E$ be a Banach lattice with $\sigma$-order continuous norm (Theorem \ref{theorem:banach_lattice_sigma_order_c_norm_mediated_and_splitting}), or $E$ is the space mentioned in \ref{observation:countable_sum_property}(II)). 
Then $\Gamma_V$ is stable and $\Gamma_{VL}=\Gamma_{LV}$, $\varphi_{VL}=\varphi_{LV}$. 
\end{theorem}
\begin{proof}
This is a consequence of  
Theorem \ref{theorem:Gamma_LV_is_subset_of_Gamma_VL} 
and 
Theorem \ref{theorem:conditions_for_VL_subset_LV}. 
\end{proof}

For a Riesz space $F$ and a Riesz subspace $\Gamma$ of $F^X$ we will now investigate under which conditions on $\varphi(\Gamma)$, $\varphi_L(\Gamma_L)$ and $\varphi_V(\Gamma_V)$ the spaces $\Gamma_{LV}$ and $\Gamma_{VL}$ are Riesz subspaces of $F^X$.

\begin{theorem}
\label{theorem:extensions_are_Riesz_spaces_if_E_is_sigma-R_complete}
Suppose $F$ is a Riesz space and $\Gamma$ is a Riesz subspace of $F^X$.
If $\varphi(\Gamma)$ is splitting in $E$ and $\varphi_L(\Gamma_L)$ is mediated in $E$, 
then $\Gamma_{LV}$ is a Riesz subspace of $F^X$. 
If $\varphi(\Gamma)$ is mediated in $E$ and $\varphi_V(\Gamma_V)$ is splitting in $E$,
then $\Gamma_{VL}$ is a Riesz subspace of $F^X$. \\
In particular, if $E$ is mediated  (and thus splitting), then both $\Gamma_{LV}$ and $\Gamma_{VL}$ are Riesz subspaces of $F^X$. 
\end{theorem}
\begin{proof}
Note first that if  $\varphi(\Gamma)$ is mediated in $E$, then $\Gamma_V$ is stable by Theorem \ref{theorem:varphi_of_Gamma_R-closed_then_L_and_V_stable}(b). 
For a proof, combine
Theorem \ref{theorem:splitting_E_implies_Gamma_L_Riesz}
and Corollary \ref{cor:mediated_E_implies_Gamma_V_Riesz}. 
\end{proof}

The next example illustrates that $\Gamma_{LV}$ is not always included in $\Gamma_{VL}$ (given that $\Gamma_V$ is stable) even if $E$ and $F$ are Riesz spaces and $\Gamma, \Gamma_{LV}, \Gamma_{VL}$ Riesz subspaces of $F^X$. 

\begin{example} 
\label{example:gamma_vl_subsetneq_gamma_lv}
\textbf{[$\Gamma_{VL} \subsetneq \Gamma_{LV} = \Gamma_{VLV}$]} \\
For an element $b= (\beta_1,\beta_2,\dots)$ of $\R^\N$ we write $b = \sum_{n\in\N} \beta_n e_n$. \\
Consider $X=\{0,1,2,\dots\}$ and $\cI = \cP(X)$. 
Let $E= c$, $F= \R^\N$, $\Omega = F^X$. 
We view the elements of $\Omega$ as sequences $(a,b_1,b_2,\dots)$ with $a,b_1,b_2,\dots \in \R^\N$. 

Define sets $\Gamma \subset \Theta \subset \Omega$ and a map $\Phi : \Theta \rightarrow \R^\N$ by 
\begin{align}
& \Theta = \{ (a,\beta_1e_1,\beta_2 e_2, \dots) : a\in c, \beta_1,\beta_2,\dots \in \R\}, \\
& \Phi(a,\beta_1e_1,\beta_2e_2,\dots) = a + \sum_{n\in\N} \beta_n e_n \qquad (a\in c, \ \beta_1,\beta_2,\dots \in \R), \\
& \Gamma = \{ (a,  \beta_1 e_1, \beta_2 e_2, \dots ) : a\in c, \  (\beta_1,\beta_2,\dots) \in c_{00} \}. 
\end{align}
Then $\Phi(\Gamma) = c = E$; let $\varphi = \Phi|_{\Gamma}$. 
From the definition it is easy to see that $\Gamma$ is stable and $\varphi$ is laterally extendable. 
We leave it to the reader to verify that $\Gamma_V= \Gamma$, 
\begin{align}
\Gamma_L = \{ (a,\ \beta_1e_1, \beta_2 e_2, \dots ): a\in c, \  (\beta_1,\beta_2,\dots ) \in c\}
\end{align}
and $\varphi_L= \Phi$ on $\Gamma_L$. 

It follows that $\Gamma_V$ is stable and $\Gamma_{VL}= \Gamma_L \subset \Gamma_{LV} = \Gamma_{VLV}$. 
We prove $\Gamma_{VL} \ne \Gamma_{LV}$. To this end, define $h\in \Omega$ by 
\begin{align}
\begin{cases}
h(n)= (-1)^n e_n \qquad \qquad (n=1,2,\dots), \\
h(0)= - \sum_{n\in\N} h(n) = - \sum_{n\in\N} (-1)^n e_n. 
\end{cases}
\end{align}
As $h(0) \notin c$ we have $h\1_{\{0\}} \notin \Gamma$; in particular, $h$ is not partially in $\Gamma$, so $h\notin \Gamma_L = \Gamma_{VL}$. It remains to prove $h\in \Gamma_{LV}$. 

For $k\in \N$, define $\tau_k, \sigma_k : X \rightarrow \R^\N$: 
\begin{align}
& \begin{cases}
\tau_k(0)= - \sum_{n=1}^k (-1)^n e_n + \sum_{n=k+1}^\infty e_n, \\
\tau_k(n)= h(n) = (-1)^n e_n  & (n=1,\dots,k), \\
\tau_k(n)= e_n 			&(n=k+1,k+2,\dots), 
\end{cases} \\
& \begin{cases}
\sigma_k(0)= - \sum_{n=1}^k (-1)^n e_n - \sum_{n=k+1}^\infty e_n,  \\
\sigma_k(n)= h(n) = (-1)^n e_n    & (n=1,\dots,k), \\
\sigma_k(n)= - e_n 			 & (n=k+1,k+2,\dots).
\end{cases}
\end{align}
Then $\tau_k, \sigma_k\in \Gamma_L$, $\tau_k \ge h \ge \sigma_k$, $\varphi_L(\tau_k) = \Phi(\tau_k) = 2 \sum_{n>k} e_n$, $\varphi_L(\sigma_k) = - 2 \sum_{n>k} e_n$, so 
$\infk \varphi_L(\tau_k) = \supk \varphi_L(\sigma_k) =0$, and $h\in \Gamma_{LV}$. 
\end{example}

The next example illustrates that $\Gamma_{VL}$ is not always included in $\Gamma_{LV}$; it provides an example of an $f\in \Gamma_{VL}^+$ for which there exist no $\rho\in \Gamma_L^+$ with $ f \le \rho$ (see Theorem \ref{theorem:extended_functions_inbetween_Gamma_functions}\emph{\ref{item:f_in_LV_iff_in_VL_and_inclosed_by_L}}).

\begin{example}
 \label{example:gamma_LV_subsetneq_Gamma_VL}
 \textbf{[$\Gamma_{LV} \subsetneq \Gamma_{VL}$]} \\
Let $E= C[0,1]$ and let $D \subset C[0,1]$ be the set of polynomials of degree $\le 2$. 
The set $D$ is order 
dense\footnote{\label{footnote:order_dense}A subspace $D$ of a partially ordered vector space $E$ is called \emph{order dense} 
in $E$ if $x= \sup \{ d\in D: d\le x\}$ (and thus $x= \inf\{d\in D: d\ge x\}$) for all $x\in E$.  
}
in $C[0,1]$ (see \cite[Example 4.4]{vGKa08}). Hence, for all $f\in E$ there exist $(g_n)_{n\in\N},(h_n)_{n\in\N}$ in $D$ with $f= \infn g_n = \supn h_n$. 
Therefore $E$ is the vertical extension of $D$ with respect to the inclusion map $D\rightarrow E$. 

Take $X = \N$, $\cI = \cP(\N), F=E=C[0,1], \Gamma = c_{00}[D] \subset F^\N =E^\N$ and let $\varphi: \Gamma \rightarrow E$ be given by $\varphi(f) = \sum_{n\in\N} f(n)$.
Since this situation is the same as in Example  \ref{example:elements_of_gammaLV_are_partially_in_gammaV} with $D_V=E$,
we have $\Gamma_V = c_{00}[E]$ and $\Gamma_{LV} \subset \Gamma_{VL}$.

Furthermore (see \ref{example:example_latterally_extendable_sequences_explained}) 
\begin{align}
\Gamma_L^+ & = \{ f \in (D^+)^\N: \sum_n f(n) \mbox{ exists in } E\}, \\
\Gamma_{VL}^+ & = \{ f \in (E^+)^\N: \sum_n f(n) \mbox{ exists in } E\}.
\end{align}
We construct an $f\in \Gamma_{VL}^+$ that is not in $\Gamma_{LV}$. 
For $n\in\N$ let $f_n$ be the `tent' function 
defined by 
\begin{align}
\notag &f_n(0)=0; \qquad f_n(\tfrac1n)=1; \qquad f_n(\tfrac1i)=0 \quad \mbox{ if } i\in\N, i\ne n; \\
&f_n \mbox{ is affine on the interval } [\tfrac{1}{1+i},\tfrac1i] \mbox{ for all } i \in \N. 
\end{align}
\begin{figure}[h]
\centering
\includegraphics[scale=1]{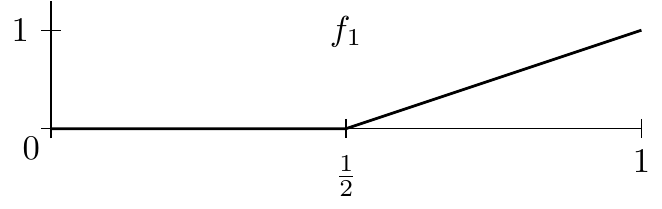}\quad \includegraphics[scale=1]{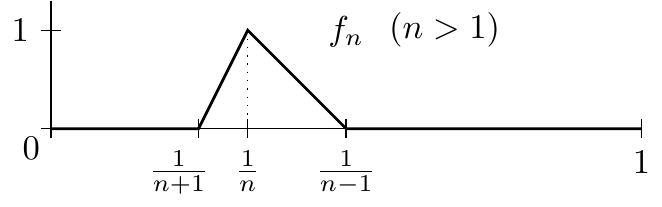}
\caption{Graph of $f_n$.}
\label{fig:picture_f_n}
\end{figure}
Then $\sum_{n=1}^\infty f_n = \1_{(0,1]}$ pointwise, so $\sum_n f_n =\1$ in $C[0,1]$.  Hence $f= (f_1,f_2,f_3,\dots)\in \Gamma_{VL}^+$. \\
We will prove that $f\notin \Gamma_{LV}$; by showing there exists no $\rho\in \Gamma_L$ for which $f \le \rho$. 

Suppose $\rho \in \Gamma_L$ and $f\le \rho$. 
Then $\rho = (\rho_1, \rho_2, \dots)$ where $\rho_1,\rho_2,\dots$ are elements of $D^+$ and $j = \sum_n \rho_n$ exists in $E= C[0,1]$. 
Let $M$ be the largest value of $j$. 
Every $\rho_n$ is a quadratic function that maps $[0,1]$ into $[0,M]$. Consequently (see the postscript)
\begin{align}
|\rho_n(x) - \rho_n(y) | \le 4 M |x-y| \qquad (x,y\in [0,1], n\in\N). 
\end{align}
In particular, $\rho_n(0) \ge \rho_n( \frac1n ) - 4M \frac1n \ge f_n(\frac1n) - 4M \frac1n= 1 - 4M \frac1n \ge \frac12$ for $n\ge 8M$. 
As $j(0) \ge \sum_{n\ge N} \rho_n(0)$ for all $N\in\N$, this is a contradiction. \\
\underline{Postscript}. Let $h : x \mapsto a x^2 +bx +c$ be a quadratic function on $[0,1]$ and $0 \le h(x) \le M$ for all $x$; we prove $|h'(x)|\le 4M$ for all $x\in [0,1]$. Since the derivative is either decreasing or increasing, we have $|h'(x)| \le \max\{|h'(0)|, |h'(1)|\}$. 
Now $h'(0)=b = 4h(\tfrac12) -h(1) -3 h(0) $ and $h'(1) = 2a+b = 3h(1) + h(0) - 4h(\tfrac12)$. Since $|h(x)-h(y)|\le M$ for all $x,y\in [0,1]$, we get the bounds $|h'(0)|\le 4M$ and $|h'(1)|\le 4M$ as desired. 
\end{example}

\begin{observation}
Observe that $\Gamma_{VL}$ in Example \ref{example:gamma_LV_subsetneq_Gamma_VL} is not stable since $(f_1,0,f_3,0,\dots)\notin \Gamma_{VL}$.
\end{observation}

\section{Embedding \texorpdfstring{$E$}{} in a (slightly) larger space}
\label{section:embeddings_in_larger_space}

\textbf{\emph{In this section $E,F,X,\cI,\Gamma,\varphi$ are as in Section \ref{section:lateral_extension}.}} 

Suppose $E^\bullet$ is another partially ordered vector space and $E \subset E^\bullet$. 
Consider $\varphi^\bullet: \Gamma \rightarrow E^\bullet$, where $\varphi^\bullet(f) = \varphi(f)$ for $f\in \Gamma$. 

Write $\Gamma_{V}^\bullet$ for the vertical extension of $\Gamma$ with respect to $\varphi^\bullet$.
If $\varphi^\bullet$ is laterally extendable, write $\Gamma^\bullet_L$ for the lateral extension of $\Gamma$ with respect to $\varphi^\bullet$, $\Gamma_{LV}^\bullet$ for the vertical extension of $\Gamma_L^\bullet$ with respect to $\varphi_L^\bullet$. 
Similarly, if $\Gamma_V^\bullet$ is stable, we introduce the notations 
$\Gamma_{VL}^\bullet$ and $\Gamma_{VLV}^\bullet$.

It is not generally the case that $\Gamma_V \subset \Gamma_V^\bullet$ or $\Gamma_L \subset \Gamma_L^\bullet$, but a natural restriction on $E^\bullet$ helps; see Theorem \ref{theorem:embed_vertical_and_lateral}.

For $E^\bullet$ we can choose to be a Dedekind complete Riesz space in which countable suprema of $E$ are preserved, in case  $E$ is integrally closed and directed (see \ref{observation:bigger_space_directed_is_fine}). In this situation, in some sense, $\Gamma_{VL}^\bullet$ is the largest extension one can obtain.

\begin{definition}
\label{def:preservation_sup_and_inf}
Let $D$ be a subspace of a partially ordered vector space $P$. 
Then we say that
\emph{countable suprema in $D$ are preserved in $P$}
if the following implication holds for all $a\in D$ and all countable $A\subset D$
\begin{align}
\label{eqn:preservation_sup}
A \mbox{ has supremum } a \mbox{ in } D \Longrightarrow
A \mbox{ has supremum } a \mbox{ in } P.
\end{align}
Note that the reverse implication holds always. 
\end{definition}

The following theorem is a natural consequence.

\begin{theorem}
\label{theorem:embed_vertical_and_lateral}
Suppose that countable suprema in $E$ are preserved in $E^\bullet$. 
Then $\varphi^\bullet$ is laterally extendable and
\begin{align}
& f\in \Gamma_V \iff f\in \Gamma_{V}^\bullet \mbox{ and } \varphi^\bullet_V(f) \in E, \\
& f\in \Gamma_L \iff f\in \Gamma_{L}^\bullet \mbox{ and } \varphi^\bullet_L(f) \in E, \\
& \varphi_V^\bullet(f) = \varphi_V(f) \mbox{ for } f\in \Gamma_V, \qquad \varphi_L^\bullet(f) = \varphi_L(f) \mbox{ for } f\in \Gamma_L, \\
\label{eqn:gamma_LV_subset_gamma_bullet_LV}
& \Gamma_{LV} \subset \Gamma_{LV}^\bullet, \qquad \varphi_{LV}^\bullet(f) = \varphi_{LV}(f) \mbox{ for } f\in \Gamma_{LV}. 
\end{align}
Suppose $\Gamma_V$ and $\Gamma_V^\bullet$ are stable. 
Then 
\begin{align}
\label{eqn:gamma_VL_subset_gamma_bullet_VL}
&\Gamma_{VL} \subset \Gamma_{VL}^\bullet, \qquad \varphi_{VL}^\bullet(f) = \varphi_{VL}(f) \mbox{ for } f\in \Gamma_{VL}, \\
\label{eqn:gamma_VLV_subset_gamma_bullet_VLV}
&\Gamma_{VLV} \subset \Gamma_{VLV}^\bullet, \qquad \varphi_{VLV}^\bullet(f) = \varphi_{VLV}(f) \mbox{ for } f\in \Gamma_{VLV}. 
\end{align}
\end{theorem}

\begin{observation}
\label{observation:bigger_space_directed_is_fine}
Under the assumptions made in \S\ref{section:lateral_extension} $\Gamma$ is directed, 
thus so are $\Gamma_L$, $\Gamma_V$ (see \ref{observation:Gamma_directed_then_Gamma_V_too}) and $\Gamma_{LV}$ (etc.). 
Hence $\varphi_V(\Gamma_V)$, $\varphi_L(\Gamma_L)$, $\varphi_{LV}(\Gamma_{LV})$ (etc.) are all subsets of $E^+ - E^+$. 
For this reason we may assume that $E$ itself is directed. 

Then under the (rather general) assumption that $E$ is also integrally closed (see Definition \ref{definition:integrally_closed}), $E$ can be embedded in a Dedekind complete Riesz space such that suprema and infima in $E$ are preserved, as we state in Theorem \ref{theorem:int_closed_povs_can_be_embedded}. 

Consequently, choosing such a Dedekind complete Riesz space for $E^\bullet$ one has the following: 
$\Gamma^\bullet_V$, $\Gamma^\bullet_{LV}$, $\Gamma^\bullet_{VL}$, $\Gamma_{VLV}^\bullet$ are stable and $\Gamma^\bullet_{LV} \subset \Gamma^\bullet_{VL}=: \overline \Gamma^\bullet$, $\overline \Gamma^\bullet_L = \overline \Gamma^\bullet_V = \overline \Gamma^\bullet$ and $\overline \varphi^\bullet_L = \overline \varphi^\bullet_V = \overline \varphi^\bullet $, where $\overline \varphi^\bullet :=\varphi_{VL}^\bullet$ (see \ref{theorem:Gamma_LV_is_subset_of_Gamma_VL}). Moreover, one has \eqref{eqn:gamma_LV_subset_gamma_bullet_LV} and if $\Gamma_V$ is stable; \eqref{eqn:gamma_VL_subset_gamma_bullet_VL} and \eqref{eqn:gamma_VLV_subset_gamma_bullet_VLV}. 
For this reason one may consider $\overline \Gamma^\bullet$ and $\overline \varphi^\bullet$  instead of $\Gamma_{LV}$ and $\varphi_{LV}$, instead of $\Gamma_{LV}^\bullet$ and $\varphi^\bullet_{LV}$ or instead of $\Gamma_{VLV}$ and $\varphi_{VLV}$, indeed $\overline \Gamma^\bullet$ contains all of the other extensions and $\overline \varphi^\bullet$ agrees with all integrals. 
\end{observation}

\begin{theorem}\cite[Chapter 4, Theorem 1.19]{Pe67} 
\label{theorem:int_closed_povs_can_be_embedded}\\
Let $E$ be an integrally closed directed partially ordered vector space. 
Then $E$ can be embedded in a Dedekind complete Riesz space $\hat E$: \\
There exists an injective linear $\gamma: E \rightarrow \hat E$ for which 
\begin{enumerate}[label=\emph{(\alph*)},align=left,leftmargin=1cm,topsep=3pt,itemsep=0pt,itemindent=-1.2em,parsep=0pt]
\item $a \ge 0 \iff \gamma(a) \ge 0$, 
\item $\gamma(E)$ is order dense in $\hat E$ (for the definition of order dense see  the  \ordinaltoname{\getrefnumber{footnote:order_dense}}
 footnote). 
\end{enumerate}
Consequently, suprema in $\gamma(E)$ are preserved in $\hat E$. 
\end{theorem}

\section{Integration for functions with values in \texorpdfstring{$\R$}{}}
\label{section:integration_in_R}

\textbf{\emph{In this section $(X,\cA,\mu)$ is a complete $\sigma$-finite measure space and $E=F=\R$.}}

We write $S$ for the vector space of simple functions from $X$ to $\R$ (see  \ref{obs:simple_functions_and_their_integral}). 
Since $\R$ is a Banach lattice with $\sigma$-order continuous norm, $S_V$ is stable 
 and $S_{LV}=S_{VL}$, $\varphi_{LV}=\varphi_{VL}$ (by Theorem \ref{theorem:Gamma_VL=LV_for_banach_lattice_with_order_cont_norm}). We write $\overline S= S_{VL}$ and $\overline \varphi = \varphi_{VL}$.

\begin{theorem}
\label{theorem:R-integration_vs_classical}
$\overline S= \cL^1(\mu)$ and $\overline \varphi(f)= \int f \D \mu$ for all $f\in \overline S$.
\end{theorem}
\begin{proof}
We prove that $S_{VL}^+ \subset \cL^1(\mu)^+ \subset S_{LV}^+$ and that $\varphi_{LV}(f) = \int f \D \mu$ for all $f\in \cL^+(\mu)$.

$S_V$ consists of the
 bounded integrable functions $f$ for which $\{x\in X: f(x) \ne 0\}$ has finite measure. 
By monotone convergence,  we have $f\in \cL^1(\mu)$ for every $f\in S_{VL}^+$.

Conversely, let $f\in \cL^1(\mu)^+$; we prove $f\in S_{LV}^+$ and $\varphi_{LV}(f) = \int f \D \mu$. 
Let $t\in (1,\infty)$. 
For $n\in \Z$, put $A_n = \{x\in X: t^n \le f(x) < t^{n+1} \}$. 
Then $(A_n)_{n\in\Z}$ forms a partition. 
Define $g: = \sum_{n\in\Z} t^n \1_{A_n}$ and $h: = tg$; then $g\le f \le h$. Since 
\begin{align}
\sum_{n\in \Z} t^n \mu(A_n) \le \sum_{n\in\Z} \int f \1_{A_n} \D \mu = \int f \D \mu,
\end{align}
we have $g\in S_L$ and $\varphi_L(g) \le \int f \D \mu$. Also, $h= tg \in S_L$, and $\varphi_L(h) - \varphi_L(g) = (t-1) \varphi_L(g) \le (t-1) \int f \D \mu$. By this and Lemma \ref{lemma:V_extension_in_terms_of_infimum} it follows that $f\in S_{LV}$ and $\varphi_{LV}(f) = \int f \D \mu$. 
\end{proof}

\section{Extensions of integrals on simple functions}
\label{section:extensions_of_simple_functions}

\textbf{\emph{In this section $E$ is a directed partially ordered vector space, $(X,\cA,\mu)$ is a complete $\sigma$-finite measure space and $\cI,S,\varphi$ are as in \ref{obs:simple_functions_and_their_integral} ($F=E$).}}

In \ref{theorem:zero_function_integrable_then_integral_is_zero}--\ref{example:extension_not_zero_but_zero_integral} for $f$ in $S_{LV}$ or $S_{VL}$ we discuss the relation between $f$ being almost everywhere equal to zero and $f$ having integral zero (i.e., either $\varphi_{LV}(f)=0$ or $\varphi_{VL}(f)=0$). 

In \ref{theorem:S_V_function_times_integrable_function} we show that under some conditions a function in $S_V$ multiplied with an integrable function with values in $\R$ is a function in $S_{LV}$.

In \ref{lemma:measurable_function_times_measure_S_L_functions}--\ref{theorem:measurable_function_times_measure_S_LV_functions_reversed} we investigate the relation between the ``$LV$''-extension on simple functions with respect to $\mu$ and $\nu$, where $\nu = h \mu$ for some measurable $h: X \rightarrow [0,\infty)$. 

In \ref{theorem:LV_and_order_continuous_functions} we discuss the relation between the ``$LV$''-extension simple functions with values in $E$ or in another partially ordered vector space $F$, when one makes the composition of a function in the extension with a $\sigma$-order continuous linear map $E\rightarrow F$. 

In \ref{theorem:continuous_function_on_product_space}--\ref{example:convolution} we will prove that 
under certain conditions on $X$ the function $x\mapsto F(x,\cdot)$ is in $S_V$ for all $F\in C(X\times T)$ and we relate that to convolution of certain finite measures with continuous functions on a topological group. 

\begin{theorem}
\label{theorem:zero_function_integrable_then_integral_is_zero}
Let $f: X \rightarrow E$ and $f=0$ a.e.. 
If $f\in S_{LV}$, then $\varphi_{LV}(f)=0$. 
If $S_V$ is stable and $f\in S_{VLV}$, then $\varphi_{VLV}(f)=0$. 
\end{theorem}
\begin{proof}
Let $B= \{x\in X: f(x) \ne 0\}$. Then $B\in \cA$ and $\mu(B)=0$. \\
(I) Assume $f\in S_V$. 
Choose $\sigma, \tau \in S$ with $\sigma \le f\le \tau$. 
Then $\sigma \1_B, \tau \1_B \in S$, $\sigma \1_B \le f \le \tau \1_B$, and $\varphi(\sigma \1_B) = \varphi (\tau \1_B) =0$. Hence $\varphi_V(f) =0$. \\
(II) Suppose $\sigma \in S_L^+$ and $(A_n)_{n\in\N}$ is a $\varphi$-partition for $\sigma$. 
Then $\sigma \1_{A_n \cap B} \in S^+$ for all $n\in\N$ and $\sum_n \varphi(\sigma \1_{A_n\cap B})=0$, i.e., $\sigma \1_B \in S_L^+$ with $\varphi_L(\sigma \1_B)=0$. In particular, if $f\in S_L$ then $\varphi_L(f)=0$. \\
(III) Assume $f\in S_{LV}$. With (II) one can repeat the argument of (I) with $S$ replaced by $S_L$ and conclude $\varphi_{LV}(f)=0$. \\
(IV) Suppose $S_V$ is stable and $f\in S_{VLV}$. One can repeat the argument in (III) with $S$ replaced by $S_V$ and conclude $\varphi_{VLV}(f)=0$. 
\end{proof}

\begin{definition}
A subset $D\subset E$ is called
\emph{order bounded}
if there are $a,b\in E$ for which $a\le D \le b$. 
\end{definition}

\begin{theorem}
\label{theorem:elements_of_extension_are_order_bounded_on_partition}
Let $f\in S_{LV}$ or (assuming $S_V$ is stable) $f\in S_{VLV}$. 
Then there exists a partition $(A_n)_{n\in\N}$ such that each set $f(A_n)$ is order bounded. 
\end{theorem}
\begin{proof}
There exists a partition $(A_n)_{n\in\N}$ such that for all $n\in\N$ there exist $h_n,g_n\in S$ for which $h_n \le f\1_{A_n} \le g_n$.
Choose $a_n,b_n \in E$ for which $a_n \le  h_n(x)$ and $g_n(x)\le b_n$ for all $x\in X$. 
Then $a_n \le f(x) \le b_n$ for $n\in\N$, $x\in A_n$. 
\end{proof}

\begin{theorem}
\label{theorem:an_almost_everywhere_zero_is_integrable_when_bounded_on_a_partition}
Let $f: X \rightarrow E$ and $f=0$ a.e.. 
Suppose there exists a partition $(A_n)_{n\in\N}$ such that for every $n\in\N$ the subset $f(A_n)$ of $E$ is order bounded. Then $f\in S_{LV}$ and if $S_V$ is stable then also $f\in \Gamma_{VL}$. 
\end{theorem}
\begin{proof}
Choose $a_1,a_2,\dots$ and $b_1,b_2,\dots$ in $E$ such that 
\begin{align}
a_n \le f(x) \le b_n \qquad (n\in\N, x\in A_n). 
\end{align}
Let $B=\{x\in X: f(x) \ne 0\}$. Then $B\in \cA$ and $\mu(B)=0$. 
Hence $g:= \sum_{n\in\N} a_n \1_{A_n\cap B}$ and $h:= \sum_{n\in\N} b_n \1_{A_n \cap B}$ are  elements of $S_L$ with $\varphi(g)=0$ and $\varphi_L(h)=0$.
As $g \le f \le h$, we get $f\in S_{LV}$ and if $S_V$ is stable also $f\in S_{VL}$. 
\end{proof}

For a real valued function $f: X \rightarrow \R$ with $f\ge 0$ and $\int f \D \mu=0$ we have $f=0$ a.e.. 
We will give an example of a $f\in S_V^+$ with $\varphi_V(f) =0$ but which is nowhere zero (Example \ref{example:extension_not_zero_but_zero_integral}). 
On the positive side, in Theorem \ref{theorem:f_equal_to_inf_sigma_n_and_sup_tau_n_for_countably_generated} 
we show that $f=0$ a.e. if $f\in S_{LV}^+$ and $\varphi_{LV}(f)=0$ provided that $E$ satisfies a certain separability condition. 

\begin{definition}  
\label{def:order_dense&countably_generated}
We call a  subset $D$ of $E^+\setminus \{0\}$ 
\emph{pervasive}\footnote{Our use of the term is similar to the one of O. van Gaans and A. Kalauch do in \cite[Definition 2.3]{vGKa08B}.}
in $E$ if for all $a\in E$ with $a>0$ there exists a $d\in D$ such that $0<d\le a$. 
We say that $E$ \emph{possesses a pervasive subset} if there exists a pervasive $D\subset E^+\setminus \{0\}$. 
\end{definition}

\begin{example} 
\label{examples:spaces_with_countable_pervasive_subsets} 
The Riesz spaces $\R^\N, \ell^\infty, c,  c_0, \ell^1$ and $c_{00}$ possess countable pervasive subsets. Indeed, in each of them the set $\{ \lambda e_n: \lambda \in \Q^+, \lambda>0, n\in\N\}$ is pervasive. \\
If $\cX$ is a completely regular topological space, then $C(\cX)$ has a countable pervasive subset if and only if $\cX$ has a countable base. (If $D\subset E^+\setminus \{0\}$ is countable and pervasive, then $\fU = \{ f^{-1}(0,\infty): f\in D\}$ is a countable base; vise versa if $\fU$ is a countable base then with choosing an $f_U$ in $C(X)^+$ for each $U\in \fU$ with $f_U=0$ on $U^c$ and $f_U(x)=1$ for some $x\in U$, the set $D= \{ \epsilon f_U : \epsilon\in \Q, \epsilon>0, U\in \fU\}$ is pervasive.) \\
$L^1(\lambda)$ and $L^\infty(\lambda)$ do not possess countable pervasive subsets, considering the Lebesgue measure space $(\R,\cM,\lambda)$. (Suppose one of them does. Then one can prove the existence of non-negligible measurable sets 
$A_1,A_2,\dots \in \cM$ such that every non-negligible measurable set contains an $A_n$, whereas $\lambda(A_n) <2^{-n}$ for all $n\in\N$. Putting $C= \R \setminus \bigcup_{n\in\N} A_n$ we have a non-negligible measurable set that contains no $A_n$: a contradiction.)
\end{example}

\begin{theorem} 
\label{theorem:f_equal_to_inf_sigma_n_and_sup_tau_n_for_countably_generated}
Let $E$ possess a countable pervasive subset $D$. 
Let $f\in S_{LV}$. 
Let $ \Lambda, \Upsilon \subset S_L$ be countable sets such that $\Lambda \le f \le \Upsilon$  and $ \sup \varphi_L(\Lambda) = \inf \varphi_L(\Upsilon)$. 
Then for almost all $x\in X$
\begin{align}
\label{eqn:almost_everywhere_f}
\sup_{g\in \Lambda} g(x) = f(x) = \inf_{h\in \Upsilon} h(x). 
\end{align}
Consequently, if $f\in S_{LV}^+$ and $ \varphi_{LV} (f) =0$, then $f=0$ a.e.. (However, see Example \ref{example:extension_not_zero_but_zero_integral}.) 
\end{theorem}
\begin{proof}
(I) 
First, as a special case (namely $f=0$), let $(\tau_n)_{n\in\N} $ be a sequence in $S_L$ with $\tau_n \ge 0$ for all $n\in\N$ and $\inf_{n\in\N} \varphi_L( \tau_n )=0$. 
We prove that $\infn \tau_n(x) = 0$ for almost all $x\in X$, by proving that $\mu(A)=0$, where $A$ is the complement of the set $\{x\in X: \infn \tau_n(x)=0\}$. 
Indeed, for this $A$ we have
\begin{align}
A=  \bigcup_{d\in D} A_d, \quad \mbox{ with }  \quad A_d = \bigcap_{n\in N} \{x\in X: d\le \tau_n(x)\}.
\end{align}
Note that for all $n\in\N$ and $d\in D$ the set $\{x\in X: d\le \tau_n(x)\}$ is measurable. 
Furthermore, for all $d\in D$ we have: 
\begin{align}
 d \mu(A_d ) = \varphi( d \1_{ A_d }) \le \varphi_L(\tau_n) \qquad (n\in\N). 
\end{align}
Hence $\mu( A_d)=0$ for all $d\in D$ and thus $\mu(A) =0$. \\
(II) Suppose that $ \Lambda, \Upsilon \subset \Gamma_L$ are countable sets such that $\Lambda \le f \le \Upsilon$, $ \sup \varphi_L(\Lambda) = \inf \varphi_L(\Upsilon)$. 
Then $\inf \varphi_L(\Upsilon - \Lambda) =0$, so by (I) $\inf_{g\in \Upsilon, h\in \Lambda} (g(x) - h(x)) = 0$ for almost all $x\in X$. 
\end{proof}

\begin{example} 
\label{example:extension_not_zero_but_zero_integral}
We  give an example of a $f\in S_V^+$ with $\varphi_V(f) =0$, where $f\ne 0$ everywhere. 
Let $([0,1), \cM, \lambda)$ be the Lebesgue measure space with underlying set $[0,1)$. 
Let $E= \ell^\infty([0,1))$ (see \S\ref{section:some_notation}).
Let $f: \R \rightarrow E^+$ be defined by $f(t) = \1_{\{t\}}$ for $t\in [0,1)$. 
Note that $f$ is not partially in $S$. 
We will show $f\in S_{V}$. 
For $n\in \N$ make $\tau_n \in S$: 
\begin{align}
\tau_n (t) = \1_{[\frac{i-1}{n},\frac{i}{n})} \qquad \mbox{ if } \ i \in \{1,\dots,n\}, t\in [\tfrac{i-1}{n},\tfrac{i}{n}). 
\end{align}
Then $\varphi(\tau_n) = \frac{1}{n} \1_{[0,1)}$ and $0 \le f \le \tau_n$ for $n\in\N$, so $f\in S_V$ and $\varphi_V(f) =0$. 
But $f(t) \ne 0$ for all $t$. 
\end{example}

\begin{theorem} 
\label{theorem:S_V_function_times_integrable_function}
Let $E$ be integrally closed and mediated. 
Let $f: X \rightarrow E$ and $g: X \rightarrow \R$. We write $gf$ for the function $x\mapsto g(x) f(x)$. Then
\begin{enumerate}[label=\emph{(\alph*)},align=left,leftmargin=1cm,topsep=3pt,itemsep=0pt,itemindent=-1.2em,parsep=0pt]
\item $f\in S_{V}$ and $g$ is bounded and measurable $\Longrightarrow$ $gf\in S_{V}$. 
\item  $f$ is partially in $S_V$ and $g$ is measurable $\Longrightarrow$ $gf$ is partially in $S_V$.
\item  $f\in S_V$ and $g\in \cL^1(\mu)$ $\Longrightarrow$ $gf\in S_{LV}$.
\item  $f\in S_{VL}$ and $g$ is bounded and measurable $\Longrightarrow$ $gf\in S_{VL}$.
\item  $f\in S_{VL}$, $f(X)$ is order bounded and $g\in \cL^1(\mu)$ $\Longrightarrow$ $gf\in S_{VL}$. 
\end{enumerate}
\end{theorem}
\begin{proof} $E$ is splitting (see \ref{observation:mediated_and_splitting}\ref{observation:part_mediated_implies_splitting_in_itself}).
 \\
(a) is a consequence of Theorem \ref{theorem:bounded_measurable_functions}(a) (see also Remark \ref{remark:bounded_measurable_are_subset_of_A_V}). \\
(b) Let $(A_n)_{n\in\N}$ be a partition such that $f\1_{A_n} \in S_V$ and $g\1_{A_n}$ is bounded for all $n\in\N$. By (a) every $gf\1_{A_n}$ lies in $S_V$. Then $gf$ is partially in $S_V$. \\
(c) Assume $f\ge 0 $ and $g \ge 0$. 
Choose (see the proof of Theorem \ref{theorem:R-integration_vs_classical}) a partition $(A_n)_{n\in\N}$ and numbers $\lambda_1,\lambda_2,\dots$ in $[0,\infty)$ with 
\begin{align}
\tau:= \sumn \lambda_n \1_{A_n} \ge g, \qquad \sumn \lambda_n \mu(A_n) <\infty.
\end{align}
Then $\tau s \in S_{L}$ for all $s\in S$. 
Choose $s\in S$ with $s\ge f$. Then $0 \le gf \le \tau s$. From Theorem \ref{theorem:extended_functions_inbetween_Gamma_functions}\emph{\ref{item:splitting_then_f_in_VL_and_LV_iff_f_inclosed_by_L}} and (b) it follows that $gf\in S_{LV}$.  \\
(d) Assume $f\ge 0$ and $0\le g \le \1$. 
Using (b), choose a partition $(A_n)_{n\in\N}$ with $f\1_{A_n}\in S_V$ and $gf\1_{A_n} \in S_V$ for all $n\in\N$. Then 
\begin{align}
& 0 \le \varphi_V(gf\1_{A_n}) \le \varphi_V(f\1_{A_n}) \qquad (n\in\N).
\end{align}
Since $\sum_n \varphi_V(f \1_{A_n})$ exists and $E$ is splitting, $\sum_n \varphi_V(gf\1_{A_n})$ exists. \\
(e) Assume $f\ge 0$ and $g \ge 0$. 
Choose $a\in E^+$ with $f(x) \le a$ for all $x\in X$. 
Choose a partition $(A_n)_{n\in\N}$ and $\lambda_1,\lambda_2,\dots \in [0,\infty)$ with 
\begin{align}
& gf\1_{A_n} \in S_V \qquad (n\in\N), \\
& g \le \sumn \lambda \1_{A_n}, \qquad \sumn \lambda_n \mu(A_n) <\infty \quad \mbox{(see the proof of Theorem \ref{theorem:R-integration_vs_classical})}.
\end{align}
Then 
\begin{align}
& g f \1_{A_n} \le \lambda_n a \1_{A_n} \qquad (n\in\N), \\
& \varphi_V(\lambda_n a \1_{A_n}) = \varphi(\lambda_n a \1_{A_n}) = \lambda_n \mu(A_n) a \qquad (n\in\N),
\end{align}
so $\sum_n \varphi_V(\lambda_n a \1_{A_n})$ exists and so does $\sum_n \varphi_V(gf \1_{A_n})$. 
\end{proof}

\begin{observation}
In Lemma \ref{lemma:measurable_function_times_measure_S_L_functions}, Theorem \ref{theorem:measurable_function_times_measure_S_LV_functions} and Theorem \ref{theorem:measurable_function_times_measure_S_LV_functions_reversed} we investigate the relation between the extensions $S_{LV}$ generated by two different measures, namely $\mu$ and $h\mu$ for a measurable function $h: X \rightarrow [0,\infty)$. \\
Note that for such a function $h$ and all $s\in (1,\infty)$ there exists a $j: X \rightarrow [0,\infty)$ that is partially in the space of simple functions $X \rightarrow [0,\infty)$, i.e., $j= \sumn \alpha_n \1_{A_n}$ for a partition $(A_n)_{n\in\N}$ and $(\alpha_n)_{n\in\N}$ in $[0,\infty)$ (or in the language of \ref{observation:order_limits_step_functions} $j$ is partially in $[\cA]$) for which $j \le h \le sj$. 
In the following (\ref{lemma:measurable_function_times_measure_S_L_functions},  \ref{theorem:measurable_function_times_measure_S_LV_functions} and  \ref{theorem:measurable_function_times_measure_S_LV_functions_reversed}) we will write $\cI^\mu$, $S^\mu$ and $\varphi^\mu$ instead of $\cI$, $S$ and $\varphi$ and, similarly for another measure $\nu$ on $(X,\cA)$, 
we write $\cI^\nu,S^\nu$ and $\varphi^\nu$ 
according to \ref{obs:simple_functions_and_their_integral} with $\nu$ instead of $\mu$. 
\end{observation}

\begin{lemma}
\label{lemma:measurable_function_times_measure_S_L_functions}
Suppose $E$ is splitting. 
Let $h: X \rightarrow [0,\infty)$ be measurable, $\nu := h \mu$.
Let $s\in (1,\infty)$ and let $j: X \rightarrow [0,\infty)$ be partially in $[\cA]$ 
and such that $j\le h \le sj$. 
Let $f\in S_L^{\nu+}$. 
Then $jf\in S^\mu_L$ and $\varphi_L^\mu(jf) \le \varphi_L^\nu(f) \le s \varphi_L^\mu(jf)$. 
\end{lemma}
\begin{proof}
Assume $(A_n)_{n\in\N}$ is a partition for $j$ and a $\varphi^\mu$-partition for $f$ (so $(A_n)_{n\in\N}$ is in $\cI^\nu \cap \cI^\mu$, i.e., $\mu(A_n), \nu(A_n)<\infty$ for all $n\in\N$). 
Choose $(\alpha_n)_{n\in\N}$ in $[0,\infty)$ and $(b_n)_{n\in\N}$ in $E^+$ such that 
\begin{align}
j = \sumn \alpha_n \1_{A_n}, \qquad f = \sumn b_n \1_{A_n}. 
\end{align}
Then $jf = \sumn \alpha_n b_n \1_{A_n}$ and thus is in $S^\mu_L$ if $\sum_n \mu(A_n) \alpha_n \beta_n$ exists in $E$. 
For each $n\in\N$
\begin{align}
0 \le \mu(A_n) \alpha_n = \int j \1_{A_n} \D \mu \le \int h \1_{A_n} \D \mu = \nu(A_n),
\end{align}
whence $0 \le \mu(A_n) \alpha_n b_n \le \nu(A_n) b_n$. 
Because $f\in S^{\nu+}_L$, $\sum_n \nu(A_n) b_n$ exists in $E$. 
Since $E$ is splitting also $\sum_n \mu(A_n) \alpha_n b_n$ exists in $E$, i.e., $jf\in S^\mu_L$.

Furthermore, $\varphi^\mu_L(jf) = \sum_n \mu(A_n) \alpha_n b_n \le \sum_n \nu(A_n) b_n = \varphi_L^\nu(f)$. 
On the other hand, we get $\mu(A_n)\alpha_n = \int j \1_{A_n} \D \mu \ge \frac{1}{s} \int h \1_{A_n} \D \mu = \frac{1}{s}\nu(A_n)$ for each $n\in\N$: it follows that $\varphi^\mu_L( jf) \ge \frac{1}{s} \varphi_L^\nu(f)$. 
\end{proof}

\begin{theorem}
\label{theorem:measurable_function_times_measure_S_LV_functions}
Let $E$ be integrally closed and splitting. 
Let $h: X \rightarrow [0,\infty)$ be measurable, $\nu := h \mu$.
\begin{enumerate}[label=\emph{(\alph*)},align=left,leftmargin=1cm,topsep=3pt,itemsep=5pt,itemindent=-1.2em,parsep=0pt]
\item $f \in S_{LV}^\nu \Longrightarrow hf \in S^\mu_{LV}, \varphi^\mu_{LV}(hf) = \varphi^\nu_{LV}(f)$, 
\item $f \in  S_{VL}^\nu \Longrightarrow hf \in S_{VL}^\mu,  \varphi_{VL}^\mu(hf) =  \varphi_{VL}^\nu(f)$.
\end{enumerate}
\end{theorem}
\begin{proof}
Since both $S^\nu_{LV}$ and $ S^\nu_{VL}$ are directed, we assume $f \ge 0$. 

(a) 
Let $f\in S_{LV}^{\nu+}$. 
For $n\in\N$ let $j_n$ be partially in $[\cA]$ and such that $j_n \le h \le (1+ \frac1n) j_n$. 
Let $\Lambda, \Upsilon \subset S_L^\nu$ be countable sets with $\Lambda \le f \le \Upsilon$ be such that $\sup \varphi^\nu_L(\Lambda) = \varphi^\nu_{LV}(f) = \inf \varphi^\nu_L(\Upsilon)$. 
Then for all $\sigma \in \Lambda$ (note that $\sigma \in S_L^{\nu+}-S_L^{\nu+}$), 
$\tau \in \Upsilon$ and $n\in\N$ we have $j_n \sigma \le hf \le (1+ \frac1n) j_n \tau$ and by Lemma \ref{lemma:measurable_function_times_measure_S_L_functions} $j_n \sigma$ and $(1+\frac1n) j_n \tau$ are in $S^\mu_L$. 
Therefore we are done if both $\inf_{n\in\N,\sigma\in \Lambda, \tau \in \Upsilon} \varphi^\mu_L((1+\frac1n)j_n \tau - j_n \sigma)=0$ and $\varphi^\mu_L(j_n \sigma)\le \varphi_{LV}^\nu(f) \le \varphi^\mu_L((1+\frac1n)j_n \tau)$ for all $n\in\N$ and all $\sigma \in \Lambda, \tau \in \Upsilon$. 
By Lemma \ref{lemma:measurable_function_times_measure_S_L_functions} applied repeatedly we have
\begin{align}
0 
\notag &\le \varphi^\mu_L((1+\tfrac1n)j_n \tau - j_n \sigma) = \varphi^\mu_L(j_n \tau - j_n \sigma) + \tfrac1n \varphi^\mu_L(j_n \tau) \\
& \le \varphi^\nu_L(\tau - \sigma) + \tfrac1n \varphi_L^\nu(\tau),
\end{align}
which has infimum $0$ since $E$ is integrally closed and $\inf_{\tau\in \Upsilon,\sigma\in \Lambda} \varphi_L^\nu(\tau - \sigma) =0$. 
On the other hand, by Lemma \ref{lemma:measurable_function_times_measure_S_L_functions},
\begin{align}
\varphi^\mu(j_n \sigma) \le \varphi^\nu_L(\sigma) \le \varphi^\nu_{LV}(f) \le \varphi^\nu_L(\tau) \le (1+ \tfrac{1}{n}) \varphi^\mu_L(j_n \tau)  \qquad (n\in\N, \sigma\in \Lambda, \tau \in \Upsilon). 
\end{align}

(b) Let $f\in  S_{VL}^{\nu +}$. 
Choose a partition $(A_n)_{n\in\N}$ with $f\1_{A_n} \in S_V^\nu$ for $n\in\N$. 
By (a), $hf\1_{A_n} \in S^\mu_{LV}$ for $n\in\N$; by Lemma \ref{lemma:gamma_LV_stable_then_every_function_of_it_is_partially_V} $h f\1_{A_n}$ is partially in $S_V^\mu$. \\
Therefore we can choose a partition $(B_n)_{n\in\N}$ with 
\begin{align}
f\1_{B_n} \in S_V^\nu, \qquad h f\1_{B_n} \in S_V^\mu \qquad (n\in\N). 
\end{align}
By (a), $\varphi_V^\nu( f\1_{B_n}) = \varphi_V^\mu (hf \1_{B_n})$ for all $n\in\N$. 
But $f\in S^{\nu+}_{VL}$, so
\begin{align}
\varphi_{VL}^\nu(f) = \sum_n \varphi_V^\nu(f\1_{B_n}) = \sum_n \varphi_V^\mu(hf \1_{B_n}).
\end{align}
Then $hf \in S_{VL}^\mu$ and $\varphi^\mu_{VL}(hf) = \varphi^\nu_{VL}(f)$. 
\end{proof}

\begin{theorem}
\label{theorem:measurable_function_times_measure_S_LV_functions_reversed}
Let $E$ be integrally closed and splitting.
Let $h: X \rightarrow [0,\infty)$ be measurable, $\nu := h \mu$, $A= \{x\in X: h(x)>0\}$.
Let $f: X \rightarrow E$ be such that $hf \in S^\mu_{LV}$. Then $f\1_{A} \in S^\nu_{LV}$. 
\end{theorem}
\begin{proof}
Define $h^*: X \rightarrow [0,\infty)$ by 
\begin{align}
h^*(x) = 
\begin{cases}
\frac{1}{h(x)} & \mbox{ if } x\in A, \\
0 						& \mbox{ if } x\notin A. 
\end{cases}
\end{align}
Then $h^*$ is measurable and $hh^*= \1_{A}$ and $\1_A = \1$ $\nu$-a.e.. 

$hf$ is in $S^\mu_L$ and thus in $S_L^{\1_A\mu}$, and since $\1_A\mu = h^* \nu$, also $hf \in S_L^{h^*\nu}$. 
By Theorem \ref{theorem:measurable_function_times_measure_S_LV_functions}, 
applied to $h^*, h^*\nu, \nu, hf$ instead of $h,\nu,\mu,f$, the function $h^*h f$ is an element of $S_{LV}^\nu$. 
But $h^* h f = \1_A f $. 
\end{proof}

In Theorem \ref{theorem:LV_and_order_continuous_functions} we show that extensions of simple functions with values in $E$ composed with a $\sigma$-order continuous linear map $E\rightarrow F$ are extensions of simple functions with values in $F$ (where $E$ and $F$ are Riesz spaces). 

\begin{theorem} 
\label{theorem:LV_and_order_continuous_functions}
Let $E$ and $F$ be Riesz spaces. 
Let $S^E$ and $\varphi^E$ be as in \ref{obs:simple_functions_and_their_integral}, and let $S^F$ and $\varphi^F$ be defined analogously. 
Let $\cL_c(E,F)$ denote the set of $\sigma$-order continuous linear functions $E\rightarrow F$ and $E_c^\sim = \cL_c(E,\R)$ (definition and notation as in Zaanen \cite[Chapter 12,\S 84]{Za83}). 
Let $f\in S^E_{LV}$. 
Then $\alpha \circ f \in S^F_{LV} $ for all $\alpha \in \cL_c(E,F)$ and 
\begin{align}
\alpha \left( \varphi^E_{LV}(f) \right) = \varphi^F_{LV}( \alpha \circ f). 
\end{align}
In particular, $\alpha \circ f$ is integrable for all $\alpha \in E^\sim_c$,  and $\alpha( \varphi^E_{LV}(f)) = \int \alpha \circ f \D \mu$. 
\end{theorem}
\begin{proof}
Suppose $\alpha \in \cL_c(E,F)^+$. 
Let $\tau\in S^{E+}_L$. Suppose $\tau= \sum_{n\in \N} a_n \1_{A_n}$ for some partition $(A_n)_{n\in\N}$ and a sequence $(a_n)_{n\in\N}$ in $E^+$. 
Then
$\alpha( \varphi^E_L( \tau)) = \alpha ( \sum_n \mu(A_n)a_n ) = \sum_n \mu(A_n) \alpha(a_n)$. 
Thus $\alpha \circ \tau$ is in $S^F_{LV}$ with $\alpha( \varphi^E_L( \tau )) = \varphi^F_L (\alpha \circ \tau )$. 
Let $ (\sigma_n)_{n\in\N}, (\tau_n)_{n\in\N}$ be sequences in $S^E_L$ with $\sigma_n \le f\le \tau_n$, $ \sigma_n \uparrow, \tau_n \downarrow$ and $\varphi^E_{LV}(f) = \supn \varphi^E_L(\sigma_n)= \infn \varphi^E_L(\tau_n)$. 
Then we have $\alpha( \varphi^E_{LV}(f)) = \supn \alpha(\varphi^E_L(\sigma_n)) = \supn \varphi^F_L( \alpha \circ \sigma_n )$ and 
$\alpha( \varphi^E_{LV}(f) ) = \infn \alpha(\varphi^E_L(\tau_n)) =\infn \varphi^F_L( \alpha \circ \tau_n )$. 
Since $\alpha \circ \sigma_n  \le \alpha \circ f \le \alpha \circ \tau_n $ for all $n\in\N$, we conclude that $\alpha \circ f \in (S^F)_{LV}$ (see Theorem \ref{theorem:R-integration_vs_classical}) with $\alpha( \varphi^E_{LV}(f))= \varphi^F_{LV}(f)$. 
\end{proof}

Theorem \ref{theorem:LV_and_order_continuous_functions} will be used in \S\ref{section:comparison_B_and_P} to compare the integrals $\varphi_{LV}$ and $\varphi_{VL}$ with the Pettis integral. 

\bigskip

Before proving Theorem \ref{theorem:integration_of_one_side_of_continuous_function} we state
(in Theorem \ref{theorem:continuous_function_on_product_space})
that there is an equivalent formulation for a function $F$ to be in $C(X\times T)$ whenever $X,T$ are topological spaces and $X$ is compact. 

\begin{theorem}\cite[Theorem 7.7.5]{Se71}
\label{theorem:continuous_function_on_product_space}
Let $X$ be a compact and let $T$ be a topological space. 
Let $F: X\times T \rightarrow \R$ be such that $F(\cdot,t) \in C(X)$ for all $t\in T$. 
Then $F\in C(X\times T)$ if and only if $t\mapsto F(\cdot,t)$ is continuous, where $C(X)$ is equipped with the supremum norm. 
Consequently, if $A\subset X$ is a compact set, then $t\mapsto \sup F(A,t)$ and $t\mapsto \inf F(A,t)$ are continuous. 
\end{theorem}

\begin{theorem}
\label{theorem:integration_of_one_side_of_continuous_function}
Let $(X,d,\mu)$ be a compact metric probability space. 
Let $T$ be a topological space and $F\in C(X \times T)$. 
The function $H: X \rightarrow C(T)$ given by $H(x)= F(x,\cdot)$ is an element of $S_V$. Furthermore, for $t\in T$, $x\mapsto F(x,t)$ is integrable and 
\begin{align}
\label{eqn:R-integral_of_continuous_function_agrees_with_normal_integral}
 \left[ \varphi_V (H) \right] (t) = \int F(x,t) \D \mu(x) \qquad (t\in T). 
\end{align}
\end{theorem}
\begin{proof}
For $k\in \N$ let $A_{k1},\dots, A_{kn_k}$ be a partition of $X$ with $\diam A_{ki} \le k^{-1}$. 
Define 
\begin{align}
 \Delta_k(t) = \sup_{x,y\in X, d(x,y)<k^{-1}} |F(x,t)- F(y,t)| \qquad (t\in T). 
\end{align}
Since $x\mapsto F(x,t)$ is uniformly continuous for all $t\in T$, $\Delta_k(t) \downarrow 0$ for all $t\in T$. 
By Theorem \ref{theorem:continuous_function_on_product_space} $t\mapsto \sup F(A_{ki},t)$ and $t\mapsto \inf F(A_{ki},t)$ are continuous for all $k\in\N$ and $i\in \{1,\dots,n_k\}$. 
For $k\in\N$ let $h_k, l_k : X \rightarrow C(T)$ be given by 
\begin{align}
 h_k(x) &= \ t\mapsto \sup F(A_{ki},t) \qquad (x\in A_{ki}),\\
\notag  l_k(x) &=  \ t\mapsto \inf F(A_{ki},t) \qquad (x\in A_{ki}).
\end{align}
Then $h_k,l_k\in S$ and $(h_k(x))(t) \ge F(x,t) \ge (l_k(x))(t)$ for all $x\in X$, $t\in T$. 
For $x\in A_{ki} \cap A_{mj}$ and $t\in T$
\begin{align}
(h_k(x) -  l_m(x)  ) (t) 
& = \sup F( A_{ki},t) - \inf F(A_{mj},t) \\
\notag & \le \sup \{ F(u,t) - F(v,t) : u,v \in A_{ki}\cup A_{mj} \} \le \Delta_{k\wedge m}(t). 
\end{align}
Let $a_k = \varphi(h_k )$ and $b_k = \varphi( l_k)$ for $k\in\N$. 
Then $0 \le a_k(t) -b_m(t)\le \Delta_{k\wedge m}(t)$ for all $k,m\in\N$ and $\inf_{k,m\in\N} a_k(t) - b_m(t) \le \infk \Delta_k(t) =0$. 
Since $a_k,b_k\in C(T)$ and $\supn b_n (t) = \infn a_n(t)$ for all $t\in T$, 
the function $t \mapsto \infn a_n(t) $ is continuous, i.e., $x\mapsto F(x,\cdot)$ is an element of $S_V$. Furthermore, we conclude that the function $x\mapsto F(x,t)$ is integrable (by Theorem \ref{theorem:R-integration_vs_classical}) and  conclude \eqref{eqn:R-integral_of_continuous_function_agrees_with_normal_integral}.
\end{proof}

\begin{example}
\label{example:convolution}
Consider a metrisable locally compact group $G$. 
Let $X \subset G$ be a compact set and $\mu$ be a finite (positive) measure on $\cB(X)$, the Borel-$\sigma$-algebra of $X$. 
Let $g\in C(G)$.
Define the convolution of $g$ and $\mu$ to be the function $g * \mu : G \rightarrow \R$ given by $g * \mu (t) =  \int g(tx^{-1}) \D \mu(x) $ for $t\in G$. 
For $x\in X$, let $L_x g \in C(G)$ be the function $t\mapsto g(tx^{-1})$. 
Then by Theorem \ref{theorem:integration_of_one_side_of_continuous_function}, the function $f : X \rightarrow C(G)$ given by $f(x) = L_x g$ is in $S_V$ and $  g * \mu = \varphi_V(f) \in C(G)$.
\end{example}

\section{Comparison with Bochner- and Pettis integral}
\label{section:comparison_B_and_P}

\textbf{\emph{We consider the situation of \S\ref{section:extensions_of_simple_functions}, with an $E$ 
that has the structure of a Banach lattice.}}
We write $\|\cdot\|$ for the norm on $E$ and $E'$ for the dual of $E$. 
Then, next to our $\varphi_{LV}$ (and other extensions) there are the Bochner and the Pettis integrals.  (We refer the reader to Hille and Phillips \cite[Section 3.7]{HiPh57} for background on both integrals.)
We denote the set of Bochner (Pettis) integrable functions from the measure space $(X,\cA,\mu)$ into the Banach lattice $E$ by $\fB$ ($\fP$) and the Bochner (Pettis) integral of an integrable function $f$ by $\fb(f)$ ($\fp(f)$).

\begin{observation}
By definition of the Bochner integral, where one also starts with defining the integral on simple functions: $S\subset \fB$ and $\varphi = \fb$ on $S$. 
Since $\fB \subset \fP$ and $\fb = \fp$ on $\fB$ we also have $S\subset \fP$ with $\varphi = \fp$ on $S$. 
\end{observation}

\begin{observation}
The following is used in this section. 
The Banach dual of $E$ is equal to the order dual, i.e., $E' = E^\sim$. Moreover, for $x,y\in E$ (see  de Jonge and van Rooij \cite[Theorem 10.2]{JoRo77IntroRiesz}) 
\begin{align}
x\le y \qquad \iff \qquad \alpha( x) \le \alpha(y) \mbox{ for all } \alpha \in E^{\sim+}. 
\end{align}
This implies that for a sequence $(y_n)_{n\in\N}$ and $x,y$ in $E$: 
\begin{align}
\label{eqn:order_dual_infimum_zero_implies_infimum_zero}
\infn \alpha(y_n)=0 \mbox{ for all }\alpha \in E^{\sim+} 
\quad &\Longrightarrow \quad \infn y_n =0.
\end{align}
\end{observation}

\begin{theorem}
 Let $f\in \fP^+$ and $f$ be partially in $S$. Then $f\in S_L^+$ and $\fp(f) = \varphi_L(f)$. 
\end{theorem}
\begin{proof}
Let $(A_n)_{n\in\N}$ be a partition for which $f_n := f\1_{A_n} \in S$. 
Then for every $\alpha \in E^{\sim+}$
\begin{align}
\label{eqn:alpha_of_pettis_integral_equal_to_sum}
\alpha( \fp(f) ) = \int \alpha \circ f \D \mu = \sumn \int \alpha \circ f_n \D \mu = \sumn \alpha( \varphi(f_n)). 
\end{align}
Hence $\infN \alpha( \fp(f) - \sum_{n=1}^N \varphi(f_n)) =0$ and thus $\fp(f) = \sum_n \varphi(f_n)$ (see \eqref{eqn:order_dual_infimum_zero_implies_infimum_zero}). 
\end{proof}

\begin{theorem}
Let $f\in \fP$. Then the following holds.
\begin{enumerate}[label=\emph{(\alph*)},align=left,leftmargin=1cm,topsep=3pt,itemsep=0pt,itemindent=-1.2em,parsep=0pt]
\item If $g\in S_{LV}$ and $f\le g$, then $\fp(f) \le \varphi_{LV}(g)$. 
\item If $S_V$ is stable, $g\in S_{VLV}$ and $f\le g$, then $\fp(f) \le \varphi_{VLV}(g)$. 
\end{enumerate}
Consequently, $\fp = \varphi_{LV}$ on $\fP \cap S_{LV}$, and $\fp= \varphi_{VLV}$ on $\fP \cap S_{VLV}$ if $S_V$ is stable. \\
The statements in (a) and (b) remain valid by replacing all ``$\le$'' by ``$\ge$''.
\end{theorem}
\begin{proof}
It will be clear that if $g\in S$ and $f\le g$, then $g\in \fP$ and hence $\fp(f) \le \fp(g) =\varphi(g)$. \\
If $g\in S_V$ and $f\le g$, then there exists an $\Upsilon \subset S$ with $g \le \Upsilon$ and $\varphi_V(g) = \inf \varphi(\Upsilon) = \inf \fp(\Upsilon) \ge \fp(f)$. \\
Let $g\in S_L$ and assume $f\le g$. 
Let $g_1,g_2\in S_L^+$ be such that $g= g_1 - g_2$. 
Let $(B_i)_{i\in\N}$ be a $\varphi$-partition for both $g_1$ and $g_2$. 
Write $A_n = \bigcup_{i=1}^n B_i$ for $n\in\N$. 
Let $\alpha \in E^{\sim+}$. 
$\alpha \circ (f\1_A) = (\alpha \circ f) \1_A$ for every $A\in \cA$, so that $\alpha \circ (f\1_A)$ is integrable. Thus, for $n\in\N$ we have 
\begin{align}
\notag \int (\alpha \circ f)\1_{A_n}  \D \mu 
& = \int \alpha \circ (f\1_{A_n}) \D \mu   \le \int \alpha \circ (g \1_{A_n}) \D \mu \\
\notag 
& =  \int \alpha \circ g_1 \1_{A_n} \D \mu  - \int \alpha \circ g_2 \1_{A_n} \D \mu \\
\notag 
& = \alpha(\varphi(g_1 \1_{A_n})) - \alpha(\varphi(g_2 \1_{A_n}))\\
& \le \alpha(\varphi(g_1 \1_{A_m})) - \alpha(\varphi(g_2 \1_{A_k})) \qquad \qquad (k,m\in\N, k<n<m).
\end{align}
Which implies that $\int (\alpha \circ f)\1_{A_n} \D \mu  + \alpha (\varphi(g_2 \1_{A_k})) \le \alpha( \varphi_L(g_1))$ as soon as $k<n$. 
By letting $n$ tend to $\infty$ (as $\int (\alpha \circ f)\1_{A_n} \D \mu \rightarrow \int \alpha \circ f \D \mu = \alpha(\fp(f))$), for each $k\in \N$ we obtain
\begin{align}
\alpha(\fp(f)) \le  \alpha(\varphi_L(g_1) - \varphi(g_2 \1_{A_k})).
\end{align}
This holds for all $\alpha \in E^{\sim+}$, so 
\begin{align}
\fp(f) \le  \varphi_L(g_1) - \varphi(g_2 \1_{A_k}) .
\end{align}
This, in tern is true for every $k$, so  $\fp(f)) \le \varphi_L(g)$. \\
We leave it to check that the preceding lines can be repeated with $S_V$, $S_L$ or $S_{VL}$ instead of $S$. 
\end{proof}

\begin{theorem}
\label{theorem:sigma_order_continuous_norm_inclusions_for_Pettis_and_Bochner}
Suppose $\|\cdot\|$ is $\sigma$-order continuous. 
Write $\overline S = S_{LV} = S_{VL}$ and $\overline \varphi = \varphi_{LV} = \varphi_{VL}$ (see Theorem \ref{theorem:Gamma_VL=LV_for_banach_lattice_with_order_cont_norm}). 
\begin{enumerate}[label=\emph{(\alph*)},align=left,leftmargin=1cm,topsep=3pt,itemsep=0pt,itemindent=-1.2em,parsep=0pt]
\item 
Then $\overline S \subset \fP$. 
Consequently, if $f$ is essentially separably valued and in $\overline S$, then $f\in \fB$. In particular, $S_L \subset \fB$. 
\item Suppose there exists an $\alpha \in E^{\sim+}_c$ with the property that if $b\in E$ and $b>0$, then $\alpha(b)>0$. 
Then $\fB_V \subset \fB$. Consequently, $\overline S \subset \fB$. 
\end{enumerate}
\end{theorem}
\begin{proof}
(a) Because $\|\cdot\|$ is $\sigma$-order continuous, $E'= E^\sim_c$. 
Therefore Theorem \ref{theorem:LV_and_order_continuous_functions} implies that $\overline S \subset \fP$. 

Note that $S_L \subset \fB$. 
Since $\fB$ is a Riesz ideal in the space of strongly measurable functions $X\rightarrow E$, an $f\in \overline S$ is an element of $\fB$ if it is essentially separably valued, 
since there are elements $\sigma,\tau \in S_L$ with $\sigma \le f \le \tau$ and $f$ is weakly measurable since $f\in \fP$. 

(b) Suppose $f\in \fB_V$ and $\sigma_n, \tau_n \in \fB$ are such that $\sigma_n \le f \le \tau_n$ for $n\in\N$, $\sigma_n \uparrow, \tau_n \downarrow$ and $\supn \fb(\sigma_n) = \fb_V(f) = \infn \fb(\tau_n)$. 
Then $\infn \int \alpha\circ  (\tau_n - \sigma_n)  \D \mu = \alpha ( \infn \fb( \tau_n - \sigma_n ) ) = 0$ and therefore $\alpha( \infn (\tau_n - \sigma_n)) = \infn \alpha \circ (\tau_n - \sigma_n)$ is integrable with integral equal to zero. Therefore $\infn (\tau_n - \sigma_n)=0$ a.e., hence $\tau_n \rightarrow f$ a.e.. Therefore $f$ is strongly measurable and thus $f\in \fB$ by (a). By (a) $S_L \subset \fB$, hence $\overline S = S_{LV} \subset \fB$. 
\end{proof}

\begin{lemma}
Let $E$ be a Banach lattice with an abstract L-norm (i.e., $\|a+b\|=\|a\|+\|b\|$ for $a,b\in E^+$).
\begin{enumerate}[label=\emph{(\alph*)},align=left,leftmargin=1cm,topsep=3pt,itemsep=0pt,itemindent=-1.2em,parsep=0pt]
\item 

Then 
\begin{align}
\| \fb(f)\| = \int \| f\| \D \mu \qquad (f\in \fB^+). 
\end{align}
\item $\fB_L= \fB$. 
\item There exist an $\alpha \in E^{\sim+}_c$ as in Theorem \ref{theorem:sigma_order_continuous_norm_inclusions_for_Pettis_and_Bochner}(b). Consequently $\fB_V = \fB$. 
\end{enumerate}
\end{lemma}
\begin{proof}
(a) It is clear that $\| \fb(f)\| = \int \| f\| \D \mu$ for $f\in S^+$, hence by limits for all $f\in \fB^+$. \\
(b) Suppose $f\in \fB_L^+$. Let $(A_n)_{n\in\N}$ be a $\fb$-partition for $f$, write $f_n = f\1_{A_n}$. 
Then $\| \sum_{n=1}^N f_n - f\| \rightarrow 0$, hence $f$ is strongly measurable. 
Moreover, since $\|\cdot\|$ is $\sigma$-order continuous 
$\|\sum_{n=1}^N \fb(f_n) - \fb_L(f)\| \rightarrow 0$, hence $\sum_{n=1}^N \|\fb(f_n) \| \rightarrow \fb_L(f)$. 
Using (a) we obtain $\int \|f\| \D \mu = \sumn \int \|f_n\| \D \mu = \sumn \| \fb(f_n)\|<\infty$, i.e., $f\in \fB$. \\
(c) Extend $\alpha : E^+ \rightarrow \R$ given by $\alpha(b)=\|b\|$ to a linear map on $E$. 
\end{proof}

\begin{examples}
\label{examples:comparing_with_B_and_P}
(I) Take $X=\N$, $\cA = \cP(\N)$, and let $\mu$ be the counting measure. 
We have $S= c_{00}[E]$; $S_V=c_{00}[E]$; all functions $\N \rightarrow E$ are partially in $S$; $\overline S:= S_{LV} = S_{VL}=S_L$ (see Theorem \ref{theorem:extended_functions_inbetween_Gamma_functions}\emph{\ref{item:f_in_LV_iff_in_VL_and_inclosed_by_L}}) and $\overline S^+$ consists precisely of the functions $f: \N \rightarrow E^+$ for which $\sum_n f(n)$ exists in the sense of the ordering. 
On the other hand, $f:  \N \rightarrow E$ is Bochner integrable if and only if 
 $\sum_{n=1}^\infty \|f(n)\|<\infty$. 

$\bullet$ If $\|\cdot\|$ is a $\sigma$-order continuous norm, then $\fB \subset \overline S$. 

$\bullet$ Moreover $\|\cdot\|$ is equivalent to an abstract L-norm 
 if and only if 
 $\fB= \overline S $ (since, if $\fB=\overline S$, the following holds: if $x_1,x_2,\dots \in E^+$ and $\sum_n x_n$ exists, then $\sumn \|x_n\|<\infty$, see Theorem \ref{theorem:sufficient_condition_for_AL-norm}). 
 
$\bullet$ For $E=c_0$ there exists an $f\in \fP$ that is not in $\overline S$. For example $f: \N \rightarrow c_0$ given by 
\begin{align}
f= (e_1,-e_1,e_2,-e_2,e_3,-e_3,\dots)
\label{eqn:pettis_integrable_function_not_in_S_L}
\end{align}
is Pettis integrable since $c_0'\cong \ell^1$ has basis $\{\delta_n:n\in\N\}$ where $\delta_n(x) = x(n)$ and $\sum_{m\in\N} \delta_n(f(m)) =0$ for all $m\in\N$. $c_0$ is $\sigma$-Dedekind complete and thus by Theorem \ref{theorem:splitting_E_implies_Gamma_L_Riesz} the set $\overline S$ is a Riesz space. However, $|f|$ is not in $\overline S$ and therefore neither $f$ is. 

$\bullet$ For $E=c$ there exists an $f\in \overline S$ that is not in $\fB$ and not in $\fP$: Consider for example $f: n \mapsto e_n$. It is an element of $\overline S$ but not of $\fB$. It is not even Pettis integrable. 
(Suppose it is, and its integral is $a$. 
Then for all $u\in c'$ we have $u(a)= \int u \circ f \D \mu = \sum_{n=1}^\infty u(f(n)) = \sum_{n=1}^\infty u(e_n)$.  
Letting $u$ be the coordinate functions, we see that $a(n)=1$ for all $n\in\N$; 
letting $u$ be $x\mapsto \limn x(n)$ we have a contradiction.)

(II) $\fB \not \subset S_{VLV}$. 
Let $(\R, \cM, \lambda)$ be the Lebesgue measure space. 
Let $E$ be the $\sigma$-Dedekind complete Riesz space $L^1(\lambda)$.
Let $g\in L^1(\lambda)$ be the equivalence class of the function that equals $t^{-\frac12}$ for $0<t\le 1$ and equals $0$ for other $t$. 
Let $L_x g(t) = g(t-x)$ for $x\in \R$. 
Then the function $f: \R \rightarrow L^1(\lambda)$ for which $f(x)=\1_{[0,1]}(x) L_x g$ is Bochner integrable ($f$ is continuous in the $\|\cdot\|_1$ norm 
(because $\|L_\epsilon g - g\|_1 = 2\sqrt{\epsilon}$ for $\epsilon>0$) 
and $\int \|f(x)\|_1 \D \lambda(x) = \int \int |g(t-x)| \D \lambda(t)  \D \lambda(x) =\|g\|_1 <\infty$) but no element of $S_{VLV}$ (by Theorem \ref{theorem:elements_of_extension_are_order_bounded_on_partition}).
\end{examples}

\section{Extensions of Bochner integrable functions}
\label{section:extensions_of_bochner}

\textbf{\emph{Consider the situation of \S\ref{section:comparison_B_and_P}.}}

As we have seen in Examples \ref{examples:comparing_with_B_and_P}, e.g., \eqref{eqn:pettis_integrable_function_not_in_S_L}, the set of Pettis integrable functions need not be stable. 
We show that $\fB$ is stable and $\fb$ is laterally extendable. Furthermore we give an example of an $f\in \fB_{LV}$ that is neither in $S_{VLV}$, nor in $\fB_L$ or $\fB_V$.  

\begin{theorem}
\label{theorem:Bochner_laterally_extendable}
$\fB$ is stable and $\fb$ is laterally extendable. 
\end{theorem}
\begin{proof}
Note that $f\1_B \in \fB$ for all $f\in \fB$ and $B\in \cA$  (since $f\1_B$ is strongly measurable and $\|f\1_B\|$ is integrable), i.e., $\fB$ is stable. 
Let $(A_n)_{n\in\N}$ be a partition in $\cA$ of $X$. 
Let $f: X \rightarrow E^+$ be a Bochner integrable function. 
Then $\int \|f\| \D \mu<\infty$ and with $B_n = A_1 \cup \cdots \cup A_n$ and Lebesgue's Dominated Convergence Theorem we obtain
\begin{align}
\left\| \fb(f - \sum_{n=1}^N f \1_{A_n}) \right\| \le \int \|f(x) - \1_{B_N}(x) f(x) \| \D \mu(x) \rightarrow 0. 
\end{align}
Thus 
\begin{align}
\fb( f )
= \limN \sum_{n=1}^N \fb( f \1_{A_n} )
=  \sum_{n} \fb( f \1_{A_n}).
\end{align}
We conclude that $\fb$ is laterally extendable. 
\end{proof}

\begin{observation}
Consider the situation of Example \ref{example:extension_not_zero_but_zero_integral}.
Since $S\subset \fB$ and $\varphi(h) = \fb( h )$ for $h\in S$: $f\in \fB_V$. 
The function $f$ is not essentially separably-valued (i.e., $f(X\setminus A)$ is not separable for all null sets $A\in \cA$), hence $f$ (and thus $g$) is not strongly measurable (see \cite[Theorem 3.5.2]{HiPh57}). 
Hence $f$ is not Bochner integrable, i.e., $f\in \fB_V$ but $f\notin \fB$.

In a similar way as has been shown in Example \ref{example:extension_not_zero_but_zero_integral}, one can show that $g: \R \rightarrow E^+$ defined by $g(t) = \1_{\{t\}}$ for $t\in \R$ is in $S_{LV}$. 
Then $g \in \fB_{LV}$ but $g \notin \fB_{V}$. 
\end{observation}

\begin{observation}
All $f\in \fB_L$ are strongly measurable. Therefore for $f\in \fB_L$ we have $f\notin \fB$ if and only if $\int \|f\| \D \mu = \infty$. 
\end{observation}

The following example illustrates that by extending the Bochner integrable functions one can obtain more than by extending the simple functions. 

\begin{example} \label{example:vertical_bochner_not_bochner}
 \textbf{[$\psi\in \fB_{V}, \psi\notin \fB$]} \\
Let $X=[2,3]$, let $\cA$ be the set of Lebesgue measurable subsets of $X$ and $\mu$ be the Lebesgue measure on $X$.
Let $M$ denote the set of equivalence classes of measurable functions $\R \rightarrow \R$. 
Let 
\begin{align}
E= \Big\{ f \in M: \sup_{x\in \R} \int_{x}^{x+1} |f| <\infty \Big\}, \quad
\|\cdot\|: E \rightarrow [0,\infty), \quad \| f\| = \sup_{x\in \R} \int_{x}^{x+1} |f|. 
\end{align}
Then $E$ equipped with the norm $\|\cdot\|$ is a Banach lattice.
$E$ is an ideal in $M$ and therefore $\sigma$-Dedekind complete (hence $S_V$ is stable; \ref{theorem:varphi_of_Gamma_R-closed_then_L_and_V_stable}).
The norm $\|\cdot\|$ is not $\sigma$-order continuous. 

For $a\in \R$, $c > 0$ define $S_{a,c}: X \rightarrow E^+$ by $S_{a,c}(x) =  \1_{(a+cx,\infty)}$. If $x,y\in X$ with $y>x$ then $\|S_{a,c}(x) - S_{a,c}(y)\| \le \| \1_{(a+cx,a+cy]}\| \le c|x-y|$, so $S_{a,c}$ is continuous and therefore strongly measurable.
Furthermore $\|S_{a,c}(x)\| = 1$ for all $x\in X$, i.e., $x\mapsto \|S_{a,c}(x)\|$ is integrable. 
Thus $S_{a,c}$ is Bochner integrable. 
For $d,e\in \R$ with $e>d$ the map $E \rightarrow \R$, $f \mapsto \int_d^e f$ is a continuous linear functional. 
Therefore 
\begin{align}
\int_{d}^e & \fb(S_{a,c})
=  \int_X \int_{d}^e (S_{a,c}(x))(t)  \D t \D x 
=
\int_{d}^e \int_X (S_{a,c}(x))(t) \D x  \D t.
\end{align}
Since this holds for all $d,e\in \R$ with $e>d$, for $t\in \R$ we have 
\begin{align}
\left(\fb(S_{a,c})\right)(t) = \int_X (S_{a,c}(x))(t) \D x = \int_2^3  \1_{(a+cx,\infty)}(t) \D x 
= \left( \tfrac{t-a}{c} \wedge 3 -2 \right) \vee 0.
\end{align}
For $k\in\N$ define  $r_k, R_k : X \rightarrow E$ by 
\begin{align}
R_k&:= S_{0,k},  \qquad 
r_k:= S_{0,k}- S_{1,k}.
\end{align}
For $x\in X$ and $k\in\N$, $r_k(x) = \1_{(kx,kx+1]}$ and $kx+1 <(k+1)x$. 
Define 
\begin{align}
\psi(x) := \1_{\bigcup_{k\in\N} (kx, kx+1]} =  \sum_{k\in\N} r_k(x), \qquad 
\sigma_n &:= \sum_{k=1}^n r_k , \qquad 
\tau_n := \sum_{k=1}^n r_k + R_{n+1}.
\end{align}
Note that $ \sigma_n \le \psi \le \tau_n $ and $\sigma_n, \tau_n \in \fB$ all for $n\in\N$. 
Since $E$ is $\sigma$-Dedekind complete and therefore mediated,
from the fact that 
\begin{align}
\infn \fb( \tau_n - \sigma_n) = \infn \fb( R_{n+1}) = 0,
\end{align}
it follows that $\psi\in \fB_V$. 
However, $\psi\notin \fB$ since $\psi$ is not essentially separably valued: \\
Let $x,y\in X$, $x<y$. We prove $\|\psi(x) - \psi(y) \| \ge 1$. 
For $k \in \N$:
\begin{align}
\notag k-1 \le \tfrac{1}{y-x} < k 
& \Longrightarrow
 \begin{cases}
1+(k-1)y \le k x, \\
1+kx < ky, 
 \end{cases} \\
& \Longrightarrow
 (kx, kx +1] \cap \bigcup_{i\in\N} (i y,i y+1] =\emptyset.
\end{align}
Hence $\| \psi(x)- \psi(y)\|\ge 1$ for all $x,y\in X$ with $x \ne y$. 

So $\psi$ is an element of $\fB_V$ but not of $\fB$ (and neither of $\fB_L$). 
\end{example}

\begin{example}  \textbf{[$f\in \fB_{LV}, f\notin \fB_L, f\notin \fB_V, f\notin S_{VLV}$]} \\
Let $(X,\cA,\mu)$ be the Lebesgue measure space $(\R,\cM,\lambda)$. Let $E$ and $\psi$ be as in Example \ref{example:vertical_bochner_not_bochner}. Define $u: \R \rightarrow E$ by 
\begin{align}
u(x) = \begin{cases}
\psi(x) & x\in [2,3], \\
0 & \text{otherwise}.
\end{cases}
\end{align}
Then $u$ is an element of $\fB_V$ and not of  $\fB_L$.
As we have seen in Examples \ref{examples:comparing_with_B_and_P}(II) there exists a $g$ in $L^1(\lambda)$ and thus in $E$ such that $v : x\mapsto \1_{[0,1]}(x) L_x g$ is an element of $\fB$ that is not an element of $S_{VLV}$. 
Furthermore $w: \R \rightarrow E$ given by $w(x) = \1_{(n,n+1]}$ for $x\in (n,n+1]$ is an element of $\fB_L$ and not of $\fB_V$. 
Therefore $f=u+v+w$ is an element of $\fB_{LV}$ (and thus of $\fB_{VL}$; see Theorem \ref{theorem:Gamma_LV_is_subset_of_Gamma_VL}) but is neither an element of $S_{VLV}$ nor of $\fB_V$ or $\fB_L$.
\end{example}

\section{Discussion}

Of course, to some extent our approach is arbitrary. 
We mention some alternatives, with comments.

\begin{observation}
The reader may have wondered why in our definition of the lateral extension the sets $A_n$ are required not only to be disjoint but also to cover $X$ (i.e., to form a partition). 
Without the covering of $X$ the definition remains perfectly meaningful, but the sum of two positive laterally integrable functions need not be laterally integrable, even in quite natural situations. 
(E.g., take $E=F=\R$ and $X=[0,1]$; let $\cI$ be the ring 
generated by the open intervals, $\Gamma$ the space of all Riemann integrable functions on $[0,1]$, and $\varphi$ the Riemann integral. If $f$ is the indicator of the Cantor set, then $\1-f$ is laterally integrable but $2\1 - f$ is not.)
\end{observation}

\begin{observation}
\label{observation:discussion_vertical_extension}
For the vertical extension we have, somewhat artificially, introduced a countability restriction leading us from $\varphi_v$ to $\varphi_V$; see Definition \ref{def:phi_V_integral}. 
In some sense, $\varphi_v$ would have served as well as $\varphi_V$. 
In order to get a non-void theory, however, we would need a much stronger (but analogous) condition than ``mediatedness'', restricting our world drastically. 
\end{observation}

\begin{observation}
\label{observation:discussion_other_extension_more_Daniel_like}
A different approach to both the vertical and the lateral extension, closer to Daniell and Bourbaki, could run as follows. 
Starting from the situation of \ref{observation:extra_assumptions_vertical_chapter}, call a function $X \rightarrow F^+$ ``integrable'' if there exist $f_1,f_2,\dots \in \Gamma^+$ such that 
\begin{align}
\begin{cases}
f_n \uparrow f \mbox{ in } F^X, \\
\supn \varphi(f_n) \mbox{ exists in } E, 
\end{cases}
\end{align}
then define the ``integral'' $\overline \varphi(f)$ of $f$ by 
\begin{align}
\overline \varphi(f) : = \supn \varphi(f_n) .
\end{align}
This definition is meaningful only if, in the above situation 
\begin{align}
g\in \Gamma^+, g \le f, \quad \Longrightarrow \quad \varphi(g) \le \supn \varphi(f_n)
\end{align}
which in a natural way leads to the requirement that $\Gamma$ be a lattice and that $\varphi$ be continuous in the following sense:
\begin{align}
\label{eqn:integral_continuity_property}
h_1, h_2,\dots \in \Gamma^+, h_n \downarrow 0 \quad \Longrightarrow \quad \varphi(h_n) \downarrow 0. 
\end{align}
These conditions lead to a sensible theory, but again we consider them as too restrictive. 
(See Example II.2.4 in the thesis of G.  Jeurnink \cite{je82} for an example of a $\Gamma$ that consists of simple functions on a measure space with values in a $C(X)$ for which \eqref{eqn:integral_continuity_property} does not hold for the standard integral on simple functions (see \ref{obs:simple_functions_and_their_integral}).) 
\end{observation}


\section*{Acknowledgements}

The authors are grateful to O. van Gaans for valuable discussions. 
W.B. van Zuijlen is supported by ERC Advanced Grant VARIS-267356. 


\bibliographystyle{plain}

\appendix

\section{Appendix}

\begin{theorem}
\label{theorem:sufficient_condition_for_AL-norm}
Let $E$ be a Banach lattice with the property 
\begin{align}
\mbox{ If } x_1,x_2,\dots \in E^+ \mbox{ and } \sum_n x_n \mbox{ exists, then } \sum_{n\in\N} \|x_n\|<\infty. 
\label{eqn:property_sufficient_for_AL_norm}
\end{align}
Then the norm $\|\cdot\|$ is equivalent to an L-norm. 
\end{theorem}

The proof uses the following lemma.

\begin{lemma}
\label{lemma:property_sufficient_for_AL_norm_implies_bound_on_sum_norms}
Let $E$ be a Banach lattice that satisfies \eqref{eqn:property_sufficient_for_AL_norm}. 
Then there exists a $C>0$ such that 
\begin{align}
x_1,x_2,\dots \in E^+, \sum_n x_n \mbox{ exists } \quad \Longrightarrow \quad \sum_{n\in\N} \|x_n\| \le C \Big\| \sum_n x_n \Big\|. 
\end{align}
\end{lemma}
\begin{proof}
Suppose not. For $i\in\N$ let $x_{i1}, x_{i2},\dots \in E^+$, $\sum_n x_{in} = b_i$ and $\sum_{n\in\N} \|x_{in}\|> 2^i \|b_i\|$ and $\|b_i\| = 2^{-i}$. 
Then $\sum_{i\in\N} \|b_i\| <\infty$, so $\sum_i b_i$ exists. 
As $\sum_i b_i = \sum_i \sum_n x_{in}$, by \eqref{eqn:property_sufficient_for_AL_norm} we get $\infty> \sum_{i\in\N} \sum_{n\in\N} \|x_{in}\| > \sum_{i\in\N} 2^i \|b_i\| = \infty$. 
\end{proof}

\begin{proof}[Proof of Theorem \ref{theorem:sufficient_condition_for_AL-norm}]
By Lemma \ref{lemma:property_sufficient_for_AL_norm_implies_bound_on_sum_norms} we can define $p: E \rightarrow [0,\infty)$,
\begin{align}
p(x)= \sup \left\{ \sum_{n\in\N} \|x_n\| \ : \ x_1,x_2,\dots \in E^+, \sum_n x_n \le |x| \right\},
\end{align}
obtaining 
$p(x)=p(|x|),
p(tx) = |t|p(x),
\|x\| \le p(x)\le C \|x\|$ for all $x\in E$, $t\in \R$ (with $C$ as in Lemma \ref{lemma:property_sufficient_for_AL_norm_implies_bound_on_sum_norms}) and $p(x) \le p(y)$ for $x,y\in E^+$ with $x\le y$. \\
Let $x,y\in E^+$; we prove $p(x +y) = p(x)+p(y)$. \\
$\bullet$ For $\epsilon>0$ choose $x_1,x_2,\dots, y_1,y_2,\dots \in E^+$, $\sum_n x_n \le x, \sum_n y_n \le y$, $\sumn \|x_n\| \ge p(x) - \epsilon$, $\sumn \|y_n\|\ge p(y) - \epsilon$. Considering the sequence $x_1,y_1,x_2,y_2,\dots$ we find $\sumn (\|x_n\|+\|y_n\|) \le p(x+y)$. Hence $p(x+y) \ge p(x) + p(y)$. \\
$\bullet$ On the other hand: Let $z_1,z_2,\dots \in E^+$, $\sum_n z_n \le x+y$; we prove $\sumn \|z_n\| \le p(x) + p(y)$. 
Define $u_n,v_n$ by 
\begin{align}
u_1 + \cdots + u_n = (z_1 + \cdots + z_n) \wedge x, \qquad v_n= z_n - u_n \qquad (n\in\N). 
\end{align}
Then $(z_1+\cdots + z_n) \wedge x - z_n = (z_1 + \cdots + z_n - z_n) \wedge (x- z_n) \le (z_1 + \cdots + z_{n-1}) \wedge x$, implying $u_n - z_n \le 0$; and $(z_1+ \cdots + z_n) \wedge x \ge (z_1+ \cdots + z_{n-1}) \wedge x$, implying $u_n \ge 0$. Thus 
\begin{align}
u_n \ge 0, v_n \ge 0 \qquad (n\in\N),
\end{align}
$\sumn \|u_n\| \le \sumn \|z_n\|<\infty$, so $\sum_n u_n$ exists; $\sum_n u_n \le x$, and $\sumn \|u_n\|\le p(x)$. 
$\sumn \|v_n\| \le \sumn \|z_n\| <\infty$, so $\sum_n v_n$ exists. For every $n\in\N$, $z_1+ \cdots + z_n \le (z_1+ \cdots +z_n + y) \wedge (x+ y)  = (z_1 + \cdots + z_n) \wedge x + y = u_1 + \cdots + u_n + y$, so $v_1+\cdots + v_n \le y$; then $\sum_n v_n \le y$ and $\sumn \|v_n\| \le p(y)$. \\
Thus $\sumn \|z_n\| \le \sumn \|u_n\| + \sumn \|v_n\| \le p(x) + p(y)$. 
\end{proof}

\end{document}